\documentclass[12pt,a4paper,reqno]{amsart}
\usepackage{amssymb}
\usepackage{latexsym}
\usepackage{amsmath}
\usepackage{mathrsfs}
\usepackage{cite}

\usepackage{calc}

\newcommand{\scal}[2]{\langle #1,#2\rangle}
\newcommand{\rr}[1]{\mathbf R^{#1}}

\newcommand{\nm}[2]{\Vert #1\Vert _{#2}}

\newcommand{\op}{\operatorname{Op}}
\newcommand{\sets}[2]{\{ \, #1\, ;\, #2\, \} }

\newcommand{\fy}{\varphi}
\newcommand{\cdo}{\, \cdot \, }
\newcommand{\supp}{\operatorname{supp}}

\newcommand{\wpr}{{\text{\footnotesize $\#$}}}

\newcommand{\eabs}[1]{\langle #1\rangle}     

\newcommand{\tp}{\operatorname{Tp}}
\newcommand{\vrum}{\vspace{0.1cm}}

\newcommand{\topo}{\tp _\fy (\omega )}

\numberwithin{equation}{section}          

\newtheorem{thm}{Theorem}
\numberwithin{thm}{section}
\newtheorem{prop}[thm]{Proposition}
\newtheorem{cor}[thm]{Corollary}
\newtheorem{lemma}[thm]{Lemma}
\newtheorem*{tom}{\rubrik}
\newcommand{\rubrik}{}

\theoremstyle{definition}

\newtheorem{defn}[thm]{Definition}

\theoremstyle{remark}
\newtheorem{rem}[thm]{Remark}

\newcommand{\rd}{\mathbf{R} ^{d}}
\newcommand{\rdd}{\mathbf{R} ^{2d}}
\newcommand{\Mmpq}{M^{p,q}_{(\omega )}}
\newcommand{\psdo}{pseudo-differential operator}
\newcommand{\tfa}{time-frequency analysis}
\newcommand{\intrd}{\int _{\rd }}

\def\inv{^{-1}}

\begin{document}

\begin{abstract}

We investigate the lifting property of modulation spaces  and
construct explicit isomorpisms between them.   
For each weight function $\omega$ and suitable window function $\fy $,  the
Toeplitz operator (or localization operator) $\tp _\fy (\omega )$ is
an isomorphism  from $M^{p,q}_{(\omega _0)}$ onto $M^{p,q}_{(\omega
_0/\omega )}$ for every $p,q \in [1,\infty ]$ and arbitrary weight
function $\omega _0$. The methods involve the pseudo-differential
calculus of Bony and Chemin and the Wiener algebra property of certain
symbol classes of pseudo-differential operators. 
\end{abstract}

\title{Isomorphism Properties of Toeplitz Operators and
  Pseudo-Differential Operators between Modulation Spaces}

\author{Karlheinz Gr\"ochenig}

\address{Faculty of Mathematics \\
University of Vienna \\
Nordbergstrasse 15 \\
A-1090 Vienna, Austria}

\email{karlheinz.groechenig@univie.ac.at}

\author{Joachim Toft}

\address{Department of Mathematics and Systems Engineering,
V{\"a}xj{\"o} University, Sweden}

\email{joachim.toft@vxu.se}


\date{}


\thanks{K.G. was supported by the  Marie-Curie Excellence Grant
  MEXT-CT-2004-517154} 
\maketitle

\maketitle

\section{Introduction}\label{sec0}

 Many families of function spaces possess a
lifting property. This fundamental property ensures that spaces of
similar type, but with respect to different weights, are
isomorphic, and  the isomorphism is usually  an explicit and
natural operator that is intrinsic in the definition of the function
spaces. We study isomorphisms of \psdo s and Toeplitz operators
between modulation spaces with different weights. 

The simplest example of the lifting property occurs for  the family of weighted $L^p$-spaces. Let
$\omega >0$ be a weight function on $\rd $  and $L^p_{(\omega )}(\rd
)$ be the weighted $L^p$-space defined by the norm
$\|f\|_{L^p_{(\omega )}} = \|f \, \omega \|_{L^p}$. Then $L^p(\rd )$ and $L^p_{(\omega )}(\rd
)$ are isomorphic, and the isomorphism is given explicitly by the
multiplication operator $f\to f \omega \inv $ from $L^p$ onto $L^p_{(\omega )}(\rd
)$. 

Perhaps the best known example of a lifting property is related to the
family of Besov spaces on $\rd $ (cf.~\cite{triebel83}). For fixed $p,q \in
(0,\infty )$ the (homogeneous) Besov spaces $B^{p,q}_s(\rd )$ with
smoothness parameter $s\in \mathbf{R} $ are all isomorphic, and the
isomorphism is given by a power of the Laplacian.  Precisely, the
operator $(-\Delta )^{-r/2}$ is an isomorphism from  $B^{p,q}_s(\rd )$
onto $B^{p,q}_{s+r}(\rd )$ (see ~\cite{triebel83}). 

In this paper, we study the lifting property for modulation spaces. Whereas
with Besov spaces the smoothness is measured with derivatives and
differences, the norm of a  modulation space measures the smoothness
of a function   by means of its phase-space distribution.  Precisely, 
let $\fy \in \mathscr S(\rr
d)\setminus 0$ be fixed. Then the \emph{short-time Fourier
transform} (STFT) of $f\in \mathscr S(\rr n)$ with respect to the
{\it window function} $\fy$ is defined as
\begin{equation}\label{stft}
V_\fy f(x,\xi)=(2\pi )^{-d/2}\int _{\mathbf R^d} f(y)\,  \overline
{\fy (y-x)}e^{-i\langle y,\xi\rangle}\, dy \, .
\end{equation}
Writing $\fy _{x,\xi }(y)=\fy (y-x)e^{i\langle y,\xi\rangle}$ 
and $V_\fy f(x,\xi)  = (f,\fy _{x,\xi })$, where  $(\cdo
 ,\cdo )$ denotes the scalar product on $L^2(\mathbf R^d)$,  the
 definition of the STFT can be extended  to a continuous map from
 $\mathscr S'(\mathbf R^d)$ to $\mathscr S'(\mathbf R^{2d})\cap
 C^\infty(\mathbf R^{2d})$.

\par

 For $1\leq p,q\leq \infty $ and a non-negative weight function $\omega
>0$ on $\rr{2d}$, the {\it modulation space}
$M^{p,q}_{(\omega)}(\mathbf R^d)$ is defined as the set of all $f\in
\mathscr S'(\mathbf R^d)$ such that
\begin{equation}\label{modnorm}
\|f\| _{M^{p,q}_{(\omega )}}\equiv \Big (\int \Big (\int |V_gf(x,\xi
)\omega (x,\xi )|^p\, dx\Big )^{q/p}\, d\xi \Big )^{1/q} < \infty ,
\end{equation}
(with obvious modifications when $p=\infty $ or $q=\infty$). 
 If $\omega \equiv 1$, we write  $M^{p,q}$ instead of
 $M^{p,q}_{(\omega )}$. 
See \cite{Gc2} for a systematic exposition of modulation spaces and
\cite{Fe6} for a survey. 

From abstract arguments it is known that for fixed $p,q$  each  $M^{p,q}_{(\omega)}(\mathbf
R^d)$   isomorphic to the unweighted modulation space
$M^{p,q}(\rr{d})$~\cite[Ch.~13]{Gc2}. However,  so far 
concrete isomorphisms are known only for few  special weight functions. For
instance, if $\omega (x,\xi ) = (1+|\xi |)^r = \langle \xi \rangle
^r$, then the operator $(1 - \Delta )^{-r/2}$ is an
isomorphism between $M^{p,q}_{(\langle \xi \rangle ^s)}$ and
$M^{p,q}_{(\langle \xi \rangle ^{s+r})}$. This fact resembles the
corresponding result for Besov spaces and was already established by
Feichtinger~\cite{Fe4}. More generally, if the weight function $\omega $
depends on one variable only, then it is not difficult to come up with
explicit isomorphisms. If $\omega (x,\xi ) = m(x)$ and $m$ is
sufficiently smooth, then the
multiplication  operator $f \to f m\inv $ establishes an isomorphism
between $M^{p,q}_{(m)}$ and $M^{p,q}$, which is in complete analogy
with the family of weighted $L^p$-spaces. If $\omega (x,\xi ) = \mu (\xi
)$, then the isomorphism is given by the corresponding Fourier
multipliers, in analogy with the example of Besov spaces.  (See~\cite{Te2,To6}.)

\par

In this paper we investigate the isomorphism property for modulation spaces
with general weights depending on both variables. For  the general
case the lifting property has been established only for few types of
weights. The difficulty   is that of   understanding  and
characterizing  the range of
certain pseudo-differential operators,  so-called Toeplitz
operators.  This is clearly a much harder problem than the one 
encountered for multiplication operators or Fourier multipliers. 

\par

The first hint about the concrete form of the isomorphism between
modulation spaces with different weights comes from the theory of the
Shubin classes~\cite{Sh}. These space are by definition the range of a
so-called  Toeplitz operator on
phase-space.  Informally, the Toeplitz operator with symbol $\omega $ and ``window''
$\fy$ is defined to be
\begin{equation}
  \label{eq:cc16}
  \tp  _\fy (\omega ) f  = \int _{\rr{2d}} \omega (x,\xi ) \, V_\fy f(x,\xi) \,
 \fy _{x,\xi }           \, dxd\xi \, .
\end{equation}
For precise definitions and boundedness results we refer to
Section~1, in particular Propositions~\ref{Tpcont1} and \ref{Tpcont2}.
Toeplitz operators  arise in pseudo-differential calculus~\cite{Fo,Le}, in the theory of quantization (Berezin
quantization~\cite{Berezin71}), and in signal
processing~\cite{Dau}   (under the name of 
time-frequency localization operators or STFT multipliers).

\par

The Shubin-Sobolev space $Q_s$  is defined   as the range of
the Toeplitz operator $\tp _\fy (\omega _s)$ with weight  $\omega _s (x,\xi
)=(1+|x|^2+|\xi |^2)^{s/2}$ and window $\fy (x) = e^{-|x|^2/2}$.
Thus $Q_s$ consists   of all $f\in \mathscr S'(\mathbf R^d)$ such that
$\tp _\fy (\omega )f \in {L}^2(\rr{d})$. It was shown recently
~\cite{BBR,BCG} that  Sobolev-Shubin space $Q_{(\omega _s )}(\mathbf R^d)$
coincides with the the modulation space  $M^{2,2}_{(\omega _s
  )}(\mathbf R^d)$. In other words, the Toeplitz operator $\tp _\fy
(\omega _s)$ provides an isomorphism between ${L}^2$ and $M^{2,2}_{(\omega _s
  )}(\mathbf R^d)$. Thus the modulation spaces of Hilbert type possess
the lifting property. 

\par

The lifting property between Sobolev-Shubin spaces  was extended  to more general
modulation spaces in \cite{BT}.  If  $\omega$ is smooth and  strictly
hypoelliptic and if $\fy $ is a  Schwartz function, 
then $\tp _{\fy } (\omega )$ is an isomorphism from 
$ M^{p,q}_{(\omega )}$ onto  $M^{p,q}$. In particular, $f\in
M^{p,q}_{(\omega )}$,  if and only if $\tp _\fy (\omega )f\in
M^{p,q}$. The key step in \cite{BT} involves the   Fredholm theory for elliptic
operators. Thus the  assumptions on $\omega$ and $\fy$ are essential
in this approach. 

\par

It was conjectured~\cite{CG05} (last remark)  that the  smoothness of the symbol and the window
are unnecessary for the lifting property. In fact, the following
observation   in~\cite{CG05} gave some plausibility to this conjecture:
if   $0< A \leq \omega (x,\xi ) \leq B$ for all $(x,\xi )\in \rdd $, 
then $\topo  $ is an isomorphism on $\Mmpq$. The key argument
requires  the Wiener algebra property  of certain symbol
classes. 

\par

In this paper we give a complete solution of the isomorphism problem
between modulation spaces for arbitrary moderate  weight functions. A weight on $\rr{2d}$ is called moderate, if there exist
constants $C,N\geq 0$ such that $\omega (X_1+X_2) \leq C \langle
X_1\rangle ^N \omega (X_2)$ for all $X_1,X_2\in \rr{2d}$. This
property is equivalent to the invariance of the modulation space
$\Mmpq$ under the time-frequency shifts $f \to f_{x,\xi }$, and is
therefore no restriction. 

Our main result can be stated as follows.

\begin{thm}
  Assume that $\omega $ is a moderate weight function and $\fy \in
  \mathscr{S}(\rr{d})$. Then the Toeplitz operator $\tp _\fy (\omega
  )$ is an isomorphism from $M^{p,q}_{(\omega _0)}$ onto $M^{p,q}_{(\omega _0/\omega
)}$ for every moderate weight $\omega _0$ and every $p,q\in [1,\infty ]$.
\end{thm}
In other words, the family of modulation spaces possesses the lifting
property. 

\par

We will establish several versions of this result. Firstly, the window
function may be chosen in certain modulation spaces that are much
larger than the Schwartz class. Secondly, the theorem holds for a  more
general family of modulation spaces that includes the $\Mmpq
$. Thirdly, we will also establish isomorphisms given by  pseudo-differential
operators rather than Toeplitz operators. 

\par

Our proofs use three different types of results: (1)  boundedness
properties of pseudo-differential operators and Toeplitz operators
between modulation spaces, (2) the 
pseudo-differential calculus of Bony and Chemin in~\cite{BC},  and (3) the
Wiener algebra property of  appropriate symbol classes
(cf. \cite{Beal,Gc3,Gc4,GR08,Sj2}).  

\par

The boundedness results serve to show that all operators involved are
well-defined between the respective modulation spaces. The application of these results is the
hairy technical part of our paper.

\par

The pseudo-differential calculus of Bony and Chemin~\cite{BC} enables
us to reduce the isomorphism problem to a problem involving only
operators of order $0$.  In
particular, \cite{BC} provides  a useful  set of isomorphisms between
modulation spaces of Hilbert type. 

\par

The Wiener algebra property states that certain symbol classes, in
particular the H\"ormander class $S^0_{0,0}$ and the generalized
Sj\"ostrand classes $M^{\infty ,1}_{(v)}$,  are preserved under
inversion. This property is crucial to ensure the correct boundedness
properties of the inverse of the Toeplitz operator $\tp _\fy (\omega
)$. 

\par

Let us mention that our theory covers only weights of polynomial
growth. It is possible to prove  general results for arbitrary
subexponential weights, but the arguments are rather different and
will be pursued elsewhere. 

\medspace

The paper is organized as follows: In Section~\ref{sec1} we collect the
prerequisites about modulation spaces, Toeplitz operators, and
pseudo-differential calculus. In Section~\ref{sec1.5} we prove a first
abstract isomorphism
theorem for \psdo s between modulation spaces. 

\par

In Section \ref{sec2} we will  prove that
$\tp _\fy (\omega )f\in 
M^{p,q}_{(\omega _0)}$ if and only if $f\in M^{p,q}_{(\omega _0\omega
)}$  for
arbitrary $\omega$ and $\omega _0$ (without any hypoelliptic
assumptions on the weights) and with $\fy$ belonging to an  appropriate
modulation space. If in addition $\omega$ is smooth, then the
isomorphism theorem  holds for a broader class of
modulation spaces.

\par

In Section \ref{sec3} we  construct explicit isomorphisms  between
weighted  modulation spaces.  We show that
for  an appropriate Gauss function $\Phi$ and arbitrary moderate  weights
$\omega _0$ and $\omega$, the Toeplitz operator $\tp _\fy (\omega
_0*\Phi )$ is an isomorphism 
from $M^{p,q}_{(\omega )}$ onto $M^{p,q}_{(\omega /\omega _0)}$. 

\vspace{ 3mm}

\textbf{Acknowledgement:} We would like to thank Karoline Johansson
for valuable discussions, especially concerning
Proposition~\ref{bonysobmod}.

\par

\section{Preliminaries}\label{sec1}

\par

In this section we recall some concepts from \tfa\ and  discuss some
basic results. For details we refer to the books~\cite{Fo,Gc2}. 

\vspace{ 3mm}

\subsection*{The Short-Time Fourier Transform. } The short-time Fourier
transforms (STFT) is defined by \eqref{stft}, 
whenever  $\fy \in \mathscr S(\rr d)\setminus
0$ and $f\in \mathscr S(\rr d)$. After writing $U$ for  the
map $F(x,y)\mapsto F(y,y-x)$ and $\mathscr F_2$ for  the partial
Fourier transform of $F(x,y)$ with respect to the $y$-variable, we
note that
\begin{equation}\label{stftdist}
V_\fy f (x,\xi )= \mathscr F_2(U(f\otimes \overline \fy )) (x,\xi ) =
\mathscr F(f\, \overline {\fy (\cdo -x})(\xi ).
\end{equation}
The Fourier transform $\mathscr F$ is the linear and
continuous map on $\mathscr S'(\rr d)$, which takes the form
$$
\mathscr Ff(\xi )=\widehat f(\xi )=(2\pi )^{-d/2}\int _{\rr
d}f(y)e^{-i\scal y\xi}\, dy,
$$
when $f\in \mathscr S(\rr d)$. We write $\check {f}(x) = f(-x) =
\mathscr{F}^2f(x)$.

\par

Each of the  operators $U$ and $\mathscr F_2$ is an  isomorphism on $\mathscr
S(\rr {2d})$ and possesses  a  unique extension to an  isomorphism on
$\mathscr S'(\rr {2d})$ and to a unitary operator on $L^2(\rr
{2d})$. If $\fy \in \mathscr S'(\rr d)\setminus 0$ and $f\in \mathscr
S'(\rr d)$, then we may define $V_\fy f$ by  $V_\fy f = \mathscr
F_2(U(f\otimes \overline \fy ))$. Since $\mathscr F_2$ and $U$ are
unitary bijections on $L^2(\rr {2d})$, it follows that $V_\fy f\in
L^2(\rr {2d})$, if and only if $f,\fy \in L^2(\rr d)$, and
$$
\nm {V_\fy f}{L^2(\rr {2d})} = \nm f{L^2(\rr d)}\nm \fy{L^2(\rr d)}.
$$

\par

The short-time Fourier transform is  similar to the  Wigner
distribution, which is defined as the tempered distribution 
$$
W_{f,g} (x,\xi ) =\mathscr F \big ( f(x+\cdo /2)\overline{g(x-\cdo
/2)}\, \big )(\xi )
$$
for $f,g \in \mathscr S'(\rr d)$.
By straight-forward computations it follows that
$$
V_g f(x,\xi ) =2^{-d}e^{-i\scal x\xi /2}W_{\check f,g}(-x/2,\xi
/2)\, .
$$
 If in addition $f,g\in L^2(\rr d)$, then $W_{f,g}$ is given by
 integral 
form
$$
W_{f,g} (x,\xi ) =(2\pi )^{-d/2} \int _{\rr d} 
f(x+y/2)\overline{g(x-y/2)}e^{-i\scal y\xi }\, dy.
$$

\par

Later on we shall use Weyl calculus, and then it is convenient to use
the symplectic Fourier transform and symplectic STFT which are defined
on tempered distributions on the phase space $\rr {2d}$. The
(standard) symplectic form on $\rr {2d}$ is defined by the formula
$\sigma (X,Y)=\scal y\xi -\scal x\eta$, when $X=(x,\xi )\in \rr {2d}$
and $Y=(y,\eta )\in \rr {2d}$. The symplectic Fourier transform
$\mathscr F_\sigma$ is the linear, continuous and bijective mapping on
$\mathscr S(\rr {2d})$, defined by
$$
\mathscr F_\sigma b(Y)=\pi ^{-d}\int _{\rr {2d}} b(Z)e^{2i\sigma
(Y,Z)}\, dZ,\qquad b\in \mathscr S(\rr {2d}).
$$
Again,  $\mathscr F_\sigma$ extends  to a
continuous bijection on $\mathscr S'(\rr {2d})$, and to unitary operator on
$L^2(\rr {2d})$. For each $\fy \in \mathscr S(\rr {2d})$ and $a\in
\mathscr S'(\rr {2d})$, the symplectic STFT is defined by $\mathcal
V_\fy a(X,Y)=\mathscr F_\sigma (a\, \overline {\fy (\cdo -X)})(Y)$. We
note that $\mathcal V_\fy a$ is given by
$$
\mathcal V_\fy a(X,Y) =\pi ^{-d}\int _{\rr {2d}} a(Z)\overline {\fy
(Z-X)} e^{2i\sigma (X,Z)}\, dZ,
$$
when $a\in \mathscr S(\rr {2d})$.

\par

\subsection*{Weight Functions.} A weight function is a non-negative
function  $\omega \in L^\infty _{\mathrm{loc}}(\rr d)\setminus \{ 0
\}$. Given two weights $\omega , v \in L^\infty _{\mathrm{loc}} (\rd
)$, $\omega$ is called \emph{$v$-moderate} if
\begin{equation}\label{moderate}
\omega (x_1+x_2)\le C\omega (x_1)v(x_2) \qquad \text{ for all } x_1,
x_2 \in \rr d \, 
\end{equation}
for some constant $C>0$  independent of
$x_1,x_2\in \rr d$. 

\par

We note that for $C$ fixed, then the smallest choice of $v$ which
fulfills \eqref{moderate} is submultiplicative in
the sense that $v(x_1 + x_2) \leq Cv(x_1) v(x_2)$ for $x_1,x_2 \in
\rd $ (cf. \cite{To5}). Throughout we will assume that the submultiplicative weights
are even and non-zero. Then \eqref{moderate} implies that $\omega(x)
>0$ for all $x\in \rd $ and so $1/\omega $ is always  well-defined. 

\par

If $v$ in \eqref{moderate} can be chosen as a
polynomial, then $\omega$ is called polynomially moderate. We let
$\mathscr P(\rr d)$ be the set of all polynomially moderate weight
functions. Furthermore, we let $\mathscr P_0(\rr d)$ be the set of all
$\omega \in \mathscr P(\rr d)\cap C^\infty (\rr d)$ such that
$(\partial ^\alpha \omega )/\omega \in L^\infty$.

\vspace{ 3mm}

\subsection*{Modulation Spaces.} We use the general definition of
modulation spaces taken from \cite{Fe6,FG1,Gc2}. 

\par

Assume that $\mathscr B$ is a Banach space of complex-valued
measurable functions on $\rr d$ and that $v \in \mathscr P(\rr d)$.
Then $\mathscr B$ is called a \emph{(translation) invariant
BF-space on $\rr d$} (with respect to $v$), if there is a constant $C$
such that the following conditions are fulfilled:
\begin{enumerate}
\item $\mathscr S(\rr d)\subseteq \mathscr
B\subseteq \mathscr S'(\rr d)$ (continuous embeddings).

\vrum

\item If $x\in \rr d$ and $f\in \mathscr B$, then $\tau _xf = f(\cdot
  -x) \in
\mathscr B$, and
\begin{equation}\label{translmultprop1}
\nm {\tau _xf}{\mathscr B}\le Cv(x)\nm {f}{\mathscr B}\text .
\end{equation}

\vrum

\item If  $f,g\in L^1_{loc}(\rr d)$,  $g\in \mathscr B$,   and $|f|
\le |g|$ almost everywhere, then $f\in \mathscr B$ and
$$
\nm f{\mathscr B}\le C\nm g{\mathscr B}\text .
$$
\end{enumerate}

\par

It follows  that if $f\in \mathscr B$ and $h\in L^\infty$, then
$f\cdot h\in \mathscr B$, and
\begin{equation}\label{multprop}
\nm {f\cdot h}{\mathscr B}\le C\nm f{\mathscr B}\nm h{L^\infty}.
\end{equation}

\par

\begin{rem}\label{newbfspaces}
Assume that $\omega _0,v,v_0\in \mathscr P(\rr d)$ are such
that $\omega$ is $v$-moderate, and assume that $\mathscr B$ is a
translation-invariant BF-space on $\rr d$ with respect to $v_0$. Also
let $\mathscr B_0$ be the Banach space which consists of all $f\in
L^1_{loc}(\rr d)$ such that $\nm f{\mathscr B_0}\equiv \nm {f\, \omega
}{\mathscr B}$ is finite. Then $\mathscr B_0$ is a translation
invariant BF-space with respect to $v_0v$.
\end{rem}

\par

\begin{defn}\label{bfspaces2}
Assume that $\mathscr B$ is a translation invariant BF-space on
$\rr {2d}$, $\omega \in \mathscr P(\rr {2d})$, and that $\fy \in
\mathscr S(\rr d)\setminus \{0\}$. Then the \emph{modulation space}
$M_{(\omega)} = M_{(\omega)}(\mathscr B)$ consists of all $f\in
\mathscr S'(\rr d)$ such that
\begin{equation}\label{modnorm2}
\nm f{M_{(\omega )}(\mathscr B)}
\equiv \nm {V_\fy f\, \omega }{\mathscr B}
\end{equation}
is finite. For $\omega =1$ we write  $M(\mathscr B)$ 
instead of $M_{(\omega)}(\mathscr B)$.
\end{defn}

\par

In Definition \ref{bfspaces2}
we may   assume without loss of generality that $\omega$ and $v$
belong to $\mathscr P_0$, because  there exists always  a weight
$\omega _0\in \mathscr P_0(\rr {2d})$ such that $C^{-1}\omega _0\le
\omega \le C\omega _0$ for some constant $C>0$ (see e.{\,}g.,
Corollary~\ref{isompseudos3}(1) or 
 ~\cite{To5}.) Therefore $M_{(\omega )}(\mathscr B)=M_{(\omega
_0)}(\mathscr B)$ with equivalent norms.

\par

Assume that $\omega \in \mathscr P(\rr {2d})$, $p,q\in
[1,\infty ]$, and let $L^{p,q}_{(\omega )}(\rr {2d})$ be the
mixed-norm space of all  $F\in
L^1_{loc}(\rr {2d})$ such that
$$
\nm F{L^{p,q}_{(\omega )}} \equiv \Big (\intrd \Big (\intrd |F(x,\xi
 )\omega (x,\xi )|^p\, dx\Big )^{q/p}\, d\xi \Big )^{1/q}<\infty \, 
$$
(with obvious modifications when $p=\infty$ or
$q=\infty$). If $\omega =1$, then we set $L^{p,q}_{(\omega )}
=L^{p,q}$. Choosing $\mathscr B = L^{p,q}_{(\omega )}$, we obtain the
standard  modulation spaces 
\begin{equation*}
M^{p,q}_{(\omega )}(\rr d) = M(L^{p,q}_{(\omega
)}(\rr {2d}))=M_{(\omega )}(L^{p,q}(\rr {2d})) \, .
\end{equation*}
(See also \eqref{modnorm}.) For convenience we use the notation
$M^p_{(\omega)}$  instead of $M^{p,p}_{(\omega
)}$. Furthermore,  for  $\omega =1$ we set
$$
M(\mathscr B) = M_{(\omega )}(\mathscr B),\quad M^{p,q} =
M^{p,q}_{(\omega )},\quad \text{and}\quad  M^{p} = M^{p}_{(\omega )} \, .
$$

\par

We use the notation $\mathcal M^{p,q}_{(\omega )}$ instead of
$M^{p,q}_{(\omega )}$, if the symplectic STFT is used in the
definition of modulation space norm. That is, if $\fy \in \mathscr
S(\rr {2d})\setminus \{ 0\}$ and $\omega \in \mathscr P(\rr {2d}\oplus
\rr {2d})$, then $\mathcal M^{p,q}_{(\omega )} (\rr {2d})$ consists of all
$a\in \mathscr S'(\rr {2d})$ such that
$$
\nm a{\mathcal M^{p,q}_{(\omega )}}\equiv \Big (\int _{\rr {2d}} \Big
(\int _{\rr {2d}}|\mathcal V_\fy a(X,Y)\omega (X,Y )|^p\, dX\Big
)^{q/p}\, dY \Big )^{1/q} < +\infty \, .
$$
The symplectic definition of modulation spaces does not yield any new
spaces. In fact, setting  $\omega
(X,Y)=\omega _0(X,-2\eta ,2y)$ for  $X\in \rr {2d}$ and $Y=(y,\eta
)\in \rr {2d}$, it follows from  the definition that $\mathcal M^{p,q}_{(\omega
)}=M^{p,q}_{(\omega _0)}$ with equivalent norms. 

\par

Finally, we note that  Definition~\ref{bfspaces2} also
include the amalgam spaces $W^{p,q}_{(\omega )}(\rr d)$ defined by the
norm
$$
\nm f{W^{p,q}_{(\omega )}} \equiv \Big (\intrd \Big (\intrd |V_\fy f
(x,\xi )\omega (x,\xi )|^q\, d\xi \Big )^{p/q}\, dx \Big
)^{1/p}<\infty  \, .
$$

\par

\subsection*{Properties of Modulations Spaces.} In the following
proposition we list some well-known properties of 
modulation spaces. See \cite[Ch.~11]{Gc2} for proofs.

\par

\begin{prop}\label{p1.4}
Let $p,q,p_j,q_j\in [1,\infty ]$ for $j=1,2$, and $\omega
,\omega _1,\omega _2,v\in \mathscr P(\rr {2d})$. Assume  that
$\omega$ is $v$-moderate and $\omega _2\le C\omega
_1$ for some constant $C>0$. Let $\mathscr B$ be a 
translation-invariant BF-space with respect to $v$. Then the
following are true:
\begin{enumerate}
\item[{\rm{(1)}}] The modulation space  $M_{(\omega )}(\mathscr B)$ is
  a Banach space with 
the norm  \eqref{modnorm2}. Let  $\psi \in M^1_{(v)}(\rr d)\setminus
\{0\}$.  Then $f\in
M_{(\omega )}(\mathscr B)$, if and only if $V_\psi f \omega \in
\mathscr{B}$. Moreover, $f\mapsto \|V_\psi f \omega \|_{
\mathscr{B}}$ is an equivalent  norm on $M_{(\omega )}(\mathscr{B})$. 

\vrum

\item[{\rm{(2)}}] If  $p_1\le p_2$ and $q_1\le q_2$,   then
$$
\mathscr S(\rr d)\hookrightarrow M^{p_1,q_1}_{(\omega _1)}(\rr
n)\hookrightarrow M^{p_2,q_2}_{(\omega _2)}(\rr d)\hookrightarrow
\mathscr S'(\rr d)\text \, .
$$

\vrum

\item[{\rm{(3)}}] The $L^2$ product $( \cdo ,\cdo )$ on $\mathscr
S$ extends to a continuous map from $M^{p,q}_{(\omega )}(\rr
n)\times M^{p'\! ,q'}_{(1/\omega )}(\rr d)$ to $\mathbf C$, and
$$
f\mapsto \sup  |( f, g) | 
$$
is an equivalent norm on $M^{p,q}_{(\omega )}(\rr d)$. Here the
supremum is taken over all $g\in \mathscr{S}(\rr d)$ such that
$\|g\|_{M^{p',q'}_{(1/\omega )}}\le 1$.

\vrum

\item[{\rm{(4)}}] If $p,q<\infty$, then $\mathscr S(\rr d)$ is dense in
$M^{p,q}_{(\omega )}(\rr d)$ and the dual space of $M^{p,q}_{(\omega
)}(\rr d)$ can be identified
with $M^{p'\! ,q'}_{(1/\omega )}(\rr d)$, through the form $(\cdo  ,\cdo
)_{L^2}$. Moreover, $\mathscr S(\rr d)$ is weakly dense in $M^{\infty
}_{(\omega )}(\rr d)$.
\end{enumerate}
\end{prop}

\par

\begin{rem}\label{extwindows}
For modulation spaces of the form $M^{p,q}_{(\omega )}$ with  fixed
$p,q\in [1,\infty ]$ the norm equivalence   in Proposition
\ref{p1.4}(1)  can be extended to a larger class of windows.  In fact,
assume  that $\omega ,v\in \mathscr P(\rr {2d})$ with  $\omega  $
being $v$-moderate and
$$
1\leq r\le \min (p,p',q,q') \, .
$$
 Let $\fy \in M^r_{(v)}(\rr d)\setminus \{0\}$. Then a tempered
 distribution  $f\in
 \mathscr{S}(\rd )$ belongs to  $M^{p,q}_{(\omega )}(\rr d)$, if and
only if $V_\fy f\in L^{p,q}_{(\omega )}(\rr {2d})$. Furthermore,
different choices of $\fy\in M^r_{(v)}(\rr d)\setminus \{0\}$ in $\nm
{V_\fy f}{L^{p,q}_{(\omega )}}$ give rise to equivalent norms.
(Cf. Proposition 3.1 in \cite{To10}.)
\end{rem}

\par

Proposition \ref{p1.4}{\,}(1) and Remark \ref{extwindows} allow us
to be rather vague concerning the particular  choice of $\fy \in
M^r_{(v)}\setminus \{0\}$ in \eqref{modnorm} and $\fy \in
M^1_{(v)}\setminus \{0\}$ in \eqref{modnorm2}. 

\par

For future references we remark that modulation spaces can be
arbitrarily close to $\mathscr S(\rr d)$ and $\mathscr S'(\rr d)$ in
the sense that
\begin{equation}\label{modtempspace}
\mathscr S(\rr d)=\bigcap _{s\in \mathbf R} M^{p,q}_{(v_s)}(\rr d),\quad
\mathscr S'(\rr d)=\bigcup _{s\in \mathbf R} M^{p,q}_{(v_s)}(\rr d),\quad
v_s(X)=\eabs X^s ,
\end{equation}
and $M^{p,q}_{(v_t)}\subseteq M^{p,q}_{(v_s)}$ as $s\le t$,  and
$X=(x,\xi )\in \rr {2d}$.

\par

\subsection*{Toeplitz Operators.}  Fix a symbol  $a\in
\mathscr S(\rr {2d})$ and a window $\fy \in \mathscr S(\rr d)$. Then the
Toeplitz operator $\tp _{\fy}(a)$ is defined by the formula
\begin{equation}\label{toeplitz}
(\tp _{\fy}(a)f_1,f_2)_{L^2(\rr d)} = (a\, V_{
\fy}f_1,V_{ \fy}f_2)_{L^2(\rr {2d})}\, ,
\end{equation}
when $f_1,f_2\in \mathscr S(\rr d)$. Obviously, $\tp _{\fy}(a)$ is
well-defined and extends uniquely to a  continuous operator from $\mathscr S'(\rr d)$ to $\mathscr
S(\rr d)$.

\par

The definition of Toeplitz operators can be extended to more general
classes of windows and symbols by using appropriate estimates for the 
short-time Fourier transforms in \eqref{toeplitz}. 

\par

We state two possible ways of extending \eqref{toeplitz}. 
The first result follows from \cite[Corollary
4.2]{CG1} and its proof, and the second result is a special case
of \cite[Theorem 3.1]{TB}. We use the notation $\mathcal
L(V_1,V_2)$ for the set of linear and continuous mappings from the
topological vector space  $V_1$ into the topological vector space
$V_2$. We also set
\begin{equation}\label{omega0t}
  \omega _{0,t}(X,Y)=v_1(2Y)^{t-1}\omega _0(X) \qquad \text{ for } X,Y
  \in \rdd \, .
\end{equation}

\par

\begin{prop}\label{Tpcont1}
Let $0\le t\le 1$, $p,q\in [1,\infty]$, and  $\omega ,\omega _0,v_1,v_0\in
\mathscr P(\rr {2d})$ be such that $v_0$ and $v_1$ are submultiplicative,
$\omega _0$ is $v_0$-moderate and $\omega $ is $v_1$-moderate. Set 
\begin{equation*}
v=v_1^tv_0 \quad \text{and}\quad \vartheta = \omega _0^{1/2} \,,  
\end{equation*}
and  let  $\omega _{0,t}$ be as in \eqref{omega0t}. Then the following
are true:

\par

\begin{enumerate}
\item The definition of $(a,\fy )\mapsto \tp _\fy (a)$ from $\mathscr
S(\rr {2d})\times \mathscr S(\rr d)$ to $\mathcal L(\mathscr S(\rr
d),\mathscr S'(\rr d))$ extends uniquely to a continuous map from
$\mathcal M^\infty _{(1/\omega _{0,t})}(\rr {2d})\times M^1_{(v)}(\rr
d)$ to $\mathcal L(\mathscr S(\rr d),\mathscr S'(\rr d))$.

\vrum

\item If $\fy \in M^{1}_{(v )}(\rr d)$ and $a\in \mathcal
M^{\infty }_{(1/\omega _{0,t})}(\rr {2d})$, then $\tp _\fy (a)$
extends uniquely to a continuous map from $M_{(\vartheta \omega
)}^{p,q}(\rr d)$ to $M_{(\omega /\vartheta )}^{p,q}(\rr d)$.
\end{enumerate}
\end{prop}

\par

\begin{prop}\label{Tpcont2}
Let $\omega ,\omega _1,\omega _2, v\in \mathscr P(\rr
{2d})$ be such that $\omega _1$ is $v$-moderate, $\omega_2$ is
$ v$-moderate and $\omega =\omega _1/\omega _2$. Then the
following are true:
\begin{enumerate}
\item The mapping $(a,\fy )\mapsto \tp _\fy (a)$ extends uniquely to
a continuous map from
$L^\infty _{(\omega )}(\rr {2d})\times M^2_{(v)}(\rr d)$ to
$\mathcal L(\mathscr S(\rr d),\mathscr S'(\rr d))$.

\vrum

\item If $a\in L^\infty _{(1/\omega )}(\rr {2d})$ and $\fy \in
M^2_{(v)}(\rr d)$, then  $\tp _\fy (a)$ extends uniquely to a
continuous operator from $M^2_{(\omega _1)}(\rr d)$ to
$M^2_{(\omega _2)}(\rr d)$.
\end{enumerate}
\end{prop}

\par

\subsection*{Pseudo-Differential Operators.} The  definition of 
Toeplitz operators can be extended even further by using
pseudo-differential calculus. Assume that $a\in \mathscr
S(\rr {2d})$, and fix $t\in \mathbf R$. Then the
pseudo-differential operator $\op _t(a)$ is the linear and
continuous operator on $\mathscr S(\rr d)$  defined by the formula
\begin{equation}\label{eq:gr}
\op _t(a)f(x) =(2\pi )^{-d}\iint a((1-t)x+ty,\xi )f(y)e^{i\scal
{x-y}\xi}\, dyd\xi .
\end{equation}
For general $a\in \mathscr S'(\rr {2d})$, the pseudo-differential
operator $\op _t(a)$ is defined as the continuous operator from
$\mathscr S(\rr d)$ to $\mathscr S'(\rr d)$ with distribution
kernel
$$
K_{t,a}(x,y) = (2\pi )^{-d/2}(\mathscr F_2^{-1}a)((1-t)x+ty,x-y).
$$
This definition makes sense, since the mappings $\mathscr F_2$ and
$F(x,y)\mapsto F((1-t)x+ty,y-x)$ are isomorphisms on $\mathscr
S'(\rr {2d})$. Furthermore, by the  Schwartz kernel theorem  the
map $a\mapsto \op _t(a)$ is a bijection from $\mathscr S'(\rr
{2d})$ onto $\mathcal L(\mathscr S(\rr d),\mathscr S'(\rr d))$.
 If $t=0$, then $\op _t(a)$ is equal to the normal (or Kohn-Nirenberg)
representation $\op (a)=a(x,D)$, and if $t=1/2$, then $\op _t(a)$
is the Weyl operator $\op ^w(a)=a^w(x,D)$ of $a$.

\par

We recall that for $s,t\in \mathbf R$ and $a_1,a_2\in \mathscr S'(\rr
{2d})$, we have
\begin{equation}\label{calculitransform}
\op _s(a_1) = \op _t(a_2) \qquad \Longleftrightarrow \qquad a_2(x,\xi
)=e^{i(t-s)\scal {D_x }{D_\xi}}a_1(x,\xi ).
\end{equation}
(Note  that the right-hand side of the latter equality makes sense,
since $e^{i(t-s)\scal {D_x }{D_\xi }}$ is a Fourier multiplier with
the smooth and bounded function $e^{i(t-s)\scal x \xi}$.)

\vspace{ 3mm}

\subsection*{Symbol Classes.} The standard symbol classes in
pseudo-differential calculus were introduced by H\"ormander and
are defined by appropriate polynomial decay conditions of the partial
derivatives
of the symbol. (See \cite{Ho1,Ho3} and the references therein.) In
particular  the class  $S^0_{0,0}$ consists of all symbols $a\in
C^\infty (\rdd )$ all of whose derivative are bounded. For $\omega \in
\mathscr P(\rr {2d})$ we also consider the weighted symbol class
$S_{(\omega )}(\rdd )$ which  consists of all $a\in C^\infty (\rr
{2d})$ such that $(\partial ^\alpha a)/\omega \in L^\infty (\rr
{2d})$.

\par

In fact, H{\"o}rmander introduced in \cite{Ho1,Ho3} a
much broader family of symbol classes with smooth symbols containing
$S_{(\omega )}$ (and the standard classes 
$S^r_{\rho ,\delta}$). Each symbol class
$S(\omega ,g)$ is parameterized by an appropriate
weight function $\omega$ and an appropriate Riemannian metric $g$
on the phase space. In \cite{Ho1,Ho3}, H{\"o}rmander also proved
several continuous results, for example that $\op ^w(a)$ is
continuous on $\mathscr S(\rr d)$, and extends uniquely to a
continuous operator on $\mathscr S'(\rr d)$ when $a\in S(\omega ,g)$.

\par

The theory was extended and improved in several ways by Bony,
Chemin and Lerner (cf. e.{\,}g. \cite{BC,BoL}). Bony and
Chemin introduced in~ \cite{BC} a family of
Hilbert spaces of Sobolev type  where each space $H(\omega ,g)$
depends on the weight $\omega$ and metric $g$. These spaces fit
the Weyl calculus well, because for appropriate
$\omega$ and
$\omega _0$ and  $a\in S(\omega ,g)$ the operator $\op ^w (a)$ is
continuous from $H(\omega _0,g)$ to $H(\omega _0/\omega ,g)$.
Furthermore, they proved that for appropriate $\omega$, there
are $a\in S(\omega ,g)$ and $b\in S(1/\omega ,g)$ such that
$$
\op ^w(a)\circ \op ^w(b) = \op ^w(b)\circ \op ^w(a) =
\operatorname{Id}_{\mathscr S'},
$$
the identity operator on $\mathscr S'(\rr d)$.

\par

The composition of Weyl operators corresponds to  the
\emph{Weyl product} (also called the twisted product in the
literature) of the involved operator symbols.  For
 $a,b\in \mathscr{S}'(\rd{2d} )$ 
the \emph{Weyl product} $a\wpr b$ is defined as the Weyl symbol of
$\op ^w(a)\circ \op ^w(b)$, whenever  this composition  makes sense
as a continuous operator from $\mathscr S(\rr d)$ to $\mathscr S'(\rr
d)$. In short, this relation can be written as
$$
c=a\wpr b \qquad \Longleftrightarrow \qquad \op ^w(c) = \op ^w(a)\circ \op ^w(b).
$$
We remark that $(a,b)\mapsto a\wpr b$ is a well-defined and continuous
mapping from $\mathscr S(\rr {2d})\times \mathscr S(\rr {2d})$ to
$\mathscr S(\rr {2d})$, since $\op ^w(a)\circ \op ^w(b)$ makes sense
as a continuous operator from $\mathscr S'(\rr d)$ to $\mathscr S(\rr
d)$.

\par

An important question in the calculus concerns of finding convenient
unique extensions of the Weyl product  to larger
 spaces. For example, for an  appropriate metric $g$  and  appropriate
weight functions  $\omega _1,\omega _2$ on $\rr {2d}$  the Weyl product on $\mathscr S$ is uniquely
extendable to a continuous mapping from $S(\omega _1,g)\times S(\omega
_2,g)$ to $S(\omega _1\omega _2,g)$ 
\cite[Thm.~18.5.4]{Ho3}. Consequently, if $a_1\in S(\omega
_1,g)$ and $a_2\in S(\omega _2,g)$, then $a_1\wpr a_2 \in S(\omega
_1\omega _2,g)$, or
\begin{equation}\label{weylprodhorm}
S(\omega _1,g)\wpr S(\omega _2,g)\subseteq S(\omega _1\omega _2,g).
\end{equation}

\par

For the case of constant metric, we will use  the following proposition.

\par

\begin{prop}\label{S(omega)prop}
Assume that $\omega _j\in \mathscr P(\rr {2d})$ for $j=0,1,2$, $s,t\in
\mathbf R$, and set $\omega _{0,N}(X,Y)=\omega _0(X)\eabs Y^{-N}$ when
$N\ge 0$ is an integer. Then the following is true:
\begin{enumerate}
\item If $a_1,a_2\in \mathscr S'(\rr {2d})$ satisfy $\op _s(a_1)=\op
_t(a_2)$, then $a_1\in S_{(\omega _0)}(\rr {2d})$ if and only if
$a_2\in S_{(\omega _0)}(\rr {2d})$.

\vrum

\item $S_{(\omega _1)}\wpr S_{(\omega _2)} \subseteq S_{(\omega
_1\omega _2)}$.

\vrum

\item $\displaystyle {S_{(\omega _0)} =\bigcap _{N\ge 0}M^{\infty
,1}_{(1/\omega _{0,N})}=\cap _{N\ge 0}\mathcal M^{\infty
,1}_{(1/\omega _{0,N})}}$.
\end{enumerate}
\end{prop}

\par

\begin{proof}
The assertion (1) is an immediate consequence of Theorem
18.5.10 in \cite{Ho3}, the assertion (2) follows from
\eqref{weylprodhorm} for the case of constant metric, and (3)  is proved in \cite{HTW} (cf.\ (2.21)
in \cite{HTW}).
\end{proof}

\par

In \tfa\ one also considers mapping properties for pseudo-differential
operators between  modulation spaces  or with
symbols in  modulation spaces. Especially we need the following two results, where the first one is a generalization of \cite[Theorem 2.1]{Ta} by Tachizawa, and the second one is a
weighted version of \cite[Theorem 14.5.2]{Gc2}. We refer to \cite[Theorem 2.2]{To9} for the proof of the first proposition and to \cite{To8} for the proof of the second one.

\par

\begin{prop}\label{p3.2}
Assume that $t\in \mathbf R$, $\omega ,\omega _0\in \mathscr P(\rr {2d})$, $a\in
S_{(\omega )}(\rr {2d})$, $t\in \mathbf R$, and that $\mathscr B$
is an invariant BF-space on $\rr {2d}$. Then
$\op _t(a)$ is continuous from $M_{(\omega _0\omega)}(\mathscr
B)$ to $M_{(\omega _0)}(\mathscr B)$, and also continuous on
$\mathscr{S}(\rr{d})$ and $\mathscr{S}'(\rr{d})$
\end{prop}

\par

\begin{prop}\label{pseudomod}
Assume that $p,q\in [1,\infty ]$,
that $\omega \in \mathscr P(\rr {2d}\oplus \rr {2d})$ and $\omega
_1,\omega _2\in \mathscr P(\rr {2d})$ satisfy
\begin{equation}\label{e5.9}
\frac {\omega _2(X-Y)}{\omega _1
(X+Y)} \le C \omega (X ,Y), \qquad \quad X,Y \in \rr{2d}, 
\end{equation}
for some constant $C$. If $a\in \mathcal M^{\infty ,1}_{(\omega )}(\rr
{2d})$, then $\op ^w(a)$ extends uniquely to a continuous map from
$M^{p,q}_{(\omega _1)}(\rr d)$ to $M^{p,q}_{(\omega _2)}(\rr d)$.
\end{prop}

\par

Finally we need the following result concerning mapping properties of
modulation spaces under the Weyl product. The result is a special case
of Theorem 0.3 in \cite{HTW}. See also ~\cite{sjo07} for related
results. 

\par

\begin{prop}\label{Weylprodmod}
Assume that $\omega _j\in \mathscr P(\rr {2d}\oplus \rr {2d})$ for $j=0,1,2$ satisfy
\begin{equation}\label{weightprodmod}
\omega _0(X,Y)\le C\omega _1(X-Y+Z,Z)\omega _2(X+Z,Y-Z),
\end{equation}
for some constant $C>0$  independent of $X,Y,Z\in \rr {2d}$. Then the
map $(a,b)\mapsto a\wpr b$ from $\mathscr S(\rr {2d})\times \mathscr
S(\rr {2d})$ to $\mathscr S(\rr {2d})$ extends uniquely to a mapping  from $\mathcal
M^{\infty ,1}_{(\omega _1)}(\rr {2d})\times \mathcal M^{\infty
  ,1}_{(\omega _2)}(\rr {2d})$ to $\mathcal M^{\infty ,1}_{(\omega
  _0)}(\rr {2d})$. 
\end{prop}

\par

In the proof of our main theorem we will need the following
consequence of Proposition~\ref{Weylprodmod}.

\par

\begin{prop} \label{corweyl}
Assume that $\omega _0, \omega , v_0, v_1 \in \mathscr{P}(\rr
{2d}\oplus \rr {2d})$, that $\omega _0$ is $v_0$-moderate and
$\omega $ is $v_1$-moderate. Set $\vartheta = \omega _0^{1/2}$, and  
\begin{eqnarray}
 \omega _1(X,Y)&=&\frac{v_0(2Y)^{1/2}v_1(2Y)}{\vartheta (X+Y) \vartheta(X-Y)} \, , \notag
\\[1ex]
\omega _2(X,Y) &=&\vartheta (X-Y)\vartheta(X+Y)v_1(2Y) \, ,
\notag
\\[1ex]
v_2(X,Y) &= & v_1(2Y)\, .
\end{eqnarray}
Then 
\begin{eqnarray}
  \label{eq:hn1}
  S_{(1/\vartheta )} \wpr \mathcal{M}^{\infty , 1} _{(\omega _1)}
  \wpr S_{(1/\vartheta )} & \subset & \mathcal{M}^{\infty , 1}
  _{(v_2)}\, , \label{ch123}
  \\[1ex]
  S_{(1/\vartheta )} \wpr \mathcal{M}^{\infty , 1} _{(v_2)}
  \wpr S_{(1/\vartheta )} & \subset & \mathcal{M}^{\infty , 1}
  _{(\omega _2)} \, . \label{ch124}
\end{eqnarray}
\end{prop}

\par

\begin{proof}
Since $S_{(1/\vartheta )} = \bigcap _{N\ge 0}\mathcal M^{\infty
,1}_{(\vartheta _N)} $ with $\vartheta _N(X,Y) = \vartheta (X)
\langle Y\rangle ^N$ (Proposition~\ref{S(omega)prop}(3)), it suffices
to argue with the symbol class $\mathcal M^{\infty , 1}
_{(\vartheta _N)}$ for some sufficiently large $N$ instead of
$S_{(1/\vartheta )}$.

\par

Introducing the intermediate weight  
$$
\omega _3(X,Y) = \frac {v_1(2Y)\vartheta (X+Y)}{\omega _0(X-Y)}.
$$
we will show that for suitable $N$
\begin{align}
\omega _3(X,Y) &\le C\omega _1(X-Y+Z,Z)\vartheta
_N(X+Z,Y-Z)\label{omegacond2} \\ 
v_1(2Y) &\le C\vartheta _N(X-Y+Z,Z)\omega
_3(X+Z,Y-Z)\label{omegacond3} \, .
\end{align}
Proposition~\ref{Weylprodmod} applied to \eqref{omegacond2} shows that
$\mathcal M^{\infty ,1} _{(\omega _1)} \wpr S_{(1/\vartheta )} \subseteq
\mathcal M^{\infty ,1} _{(\omega _3)}$, and \eqref{omegacond3} implies  that
$S_{(1/\vartheta )}  \wpr \mathcal M^{\infty ,1} _{(\omega _3)} \subseteq
\mathcal M^{\infty ,1} _{(v_2)}$,  and \eqref{ch123} now follows. 

\par

Since $\vartheta $ is $v_0 ^{1/2}$-moderate and $v_0$ grows at most
polynomially, we have $\vartheta (X-Y) \inv \leq v_0(2Z)^{1/2}
\vartheta (X-Y+2Z)\inv $ and $\vartheta (X+Y) \leq \vartheta (X+Z) \langle
Y-Z \rangle ^N$ for suitable $N\geq 0$. Using these inequalities
repeatedly in the following, a  straight-forward computation yields
\begin{align*}
\omega _3(X,Y) & = \frac {v_1(2Y)\vartheta (X+Y)^{1/2}}{\vartheta (X-Y)^2}
\\[1ex]
& \le C_1 \frac {v_0(2Z) ^{1/2}v_1(2Z)\vartheta (X+Y) \eabs
{Y-Z}^{N}}{\vartheta(X-Y+2Z)\vartheta(X-Y)}
\\[1ex]
& =C_1 \omega _1(X-Y+Z,Z)\vartheta _N(X+Z,Y-Z),
\end{align*}
for some $C_1>0$ and $N>0$.

\par

Likewise we obtain 
\begin{align*}
v_1(2Y) &= \frac {\vartheta(X-Y) v_1(2Y)\vartheta
  (X-Y)}{\vartheta(X-Y) ^2}
\\[1ex]
&\le C_1 \frac {\vartheta(X-Y) v_0(2Y)^{1/2}v_1(2Y)\vartheta (X+Y)}{\vartheta(X-Y)^2}
\\[1ex]
&\le C_2 \frac {\vartheta(X-Y+Z)\eabs
Z^Nv_0(2(Y-Z))^{1/2}v_1(2(Y-Z))\vartheta(X+Y)}{\vartheta(X-Y+2Z)^2}
\\[1ex]
& = C_2 \vartheta _N(X-Y+Z,Z)\omega _3(X+Z,Y-Z) \, .
\end{align*}

\par

The twisted convolution relation \eqref{ch124} is proved similarly. Let 
$$
\omega _4(X,Y) = \vartheta(X-Y)v_1(2Y)
$$
be the intermediate weight. Then the inequality 
\begin{align*}
\omega _4(X,Y) & = \vartheta(X-Y)v_1(2Y) \le C\omega
_0(X-Y+Z)\eabs Z^Nv_1(2(Y-Z))
\\[1ex]
&= C \vartheta _N(X-Y+Z,Z)v_2(,X+Z,Y-Z)
\end{align*}
implies that $S_{(1/\vartheta   )} \wpr \mathcal M^{\infty ,1} _{(v_2)}
\subseteq  \mathcal M^{\infty ,1} _{(\omega _4)}$. 

\par

Similarly we obtain
\begin{align*}
\omega _2(X,Y) &\le C \vartheta(X-Y)v_1(2Z)\vartheta (X+Z)\eabs {Z-Y}^N
\\[1ex]
&= C\omega _4(X-Y+Z,Z)\vartheta(X+Z)\eabs {Z-Y}^N
\\[1ex]
& = C\omega _4(X-Y+Z,Z)\vartheta _N(X+Z,Y-Z),
\end{align*}
and thus $\mathcal  M^{\infty ,1} _{(\omega _4)} \wpr S_{(1/\vartheta  )}
\subseteq \mathcal  M^{\infty ,1} _{(\omega _2)}$. 
\end{proof}

\par

\subsection*{The Wiener Algebra Property.} Proposition \ref{p3.2}
generalizes parts of Calderon-Vaillancourt's theorem. In fact, if we
let $\omega =\omega _0=1$ and $\mathscr B =L^2$, then Proposition
\ref{p3.2} asserts that $\op _t(a)$ is continuous on $L^2$ as $a\in
S^0_{0,0}$. In the same way, by letting $p=q=2$, $\omega _1=\omega
_2=1$ and $\omega (x,\xi ,\eta ,y)=v(\eta ,y)$, where $v$ is
submultiplicative, then Proposition \ref{pseudomod} shows that $\op
_t(a)$ is $L^2$-continuous when $a\in \mathcal M^{\infty ,1}_{(\omega
  )}$. 

\par

As a further crucial tool in our study of  the isomorphism property of
Toeplitz operators we need to combine these continuity results with
convenient invertibility properties. These properties are  the so-called Wiener algebra property of
certain symbol classes, and asserts that the inversion of a \psdo\
preserves the symbol class and is often refered to as the
spectral invariance of a symbol class.

\par

\begin{thm}\label{specinv}
Assume that $t\in \mathbf R$ and that $v\in \mathscr P(\rr {4d})$ is submultiplicative and
depends only on the second variable $v(X,Y)=v_0(Y)$. Then the following are true:
\begin{enumerate}
\item If  $a\in S^0_{0,0}(\rr {2d})$ and $\operatorname {Op}_t(a)$ is invertible
on $L^2(\rr d)$, then the inverse of $\operatorname {Op}_t(a)$
is equal to $\operatorname {Op}_t(b)$ for some $b \in S^0_{0,0}(\rr {2d})$. 

\vrum

\item If $a\in M^{\infty ,1}_{(v)}(\rr {2d})$ and  $\op _t(a)$ is invertible on
$L^2(\rr d)$, then the inverse of $\operatorname {Op}_t(a)$ is equal to $\op _t(b)$, for some
$b\in M^{\infty ,1}_{(v)}(\rr {2d})$.
\end{enumerate}
\end{thm}

\par

The assertion (1) in Theorem \ref{specinv} is a classical result of
Beals (cf. ~\cite{Beal}), while (2) is proved in the unweighted case
in \cite{Sj2}. We refer to \cite[Corollary 5.5]{Gc3} or \cite{Gc4} for
the general case of (2), and to \cite{GR08,GS06} for further
generalizations. We also remark that Beals' result  is an immediate consequence
of (2) because $S^0_{0,0}$ is the 
intersection of all  $M^{\infty ,1}_{(v_s)}$, with $v_s(X,Y)=\langle Y
\rangle ^s$.

\vspace{ 3mm}

\subsection*{Weyl formulation of Toeplitz operators.} We finish
this section by recalling some important relations
between Weyl operators, Wigner distributions, and Toeplitz
operators. For instance, the Weyl symbol of a Toeplitz operator
is the convolution between the Toeplitz symbol and a Wigner
distribution. More precisely,  if  $a\in \mathscr S(\rr {2d})$ and
$\fy \in \mathscr S(\rr d)$, then
\begin{equation}
\label{toeplweyl}
\tp _\fy (a) = (2\pi )^{-d/2}\op ^w(a*W_{\fy  ,\fy} ) \, .  
\end{equation}

\par

Our analysis of Toeplitz operators is
based on the pseudo-differential operator representation, given
by \eqref{toeplweyl}. Furthermore, any extension of the definition of
Toeplitz operators to cases which are not covered by Propositions
\ref{Tpcont1} and \ref{Tpcont2} are based on this representation. Here
we remark that this leads to situations were certain mapping
properties for the pseudo-differential operator representation make
sense, while similar interpretations are difficult or impossible to
make in the framework of \eqref{toeplitz} (see Remark
\ref{extensionremark} in Section \ref{sec2}). We refer to \cite{To9}
or Section \ref{sec2} for extensions of Toeplitz operators in context
of pseudo-differential operators.

\par

\section{Identifications of modulation spaces}\label{sec1.5}

\par

In this section we show that for each $\omega$ and $\mathscr B$,
there are canonical ways to identify the modulation space
$M_{(\omega )}(\mathscr B)$ with $M(\mathscr B)$ by means of
convenient bijections in the form of explicit pseudo-differential operators.
More precisely we have the following result, related to Proposition \ref{p3.2}.

\par

\begin{thm}\label{identification}
Assume that $\omega \in \mathscr P(\rr {2d})$ and $t\in \mathbf
R$. Then the following are  true:
\begin{enumerate}
\item There exist $a\in S_{(\omega )}(\rr {2d})$ and $b\in S_{(1/\omega
)}(\rr {2d})$ such that
\begin{equation}\label{abinverse}
\op _t(a)\circ \op _t(b) =\op _t(b)\circ \op _t(a) =\operatorname
{Id}_{\mathscr S'(\rr d)}.
\end{equation}
Furthermore, $\op _t(a)$ is an isomorphism from $M_{(\omega
_0)}(\mathscr B)$ onto $M_{(\omega _0/\omega )} (\mathscr B)$, for
every $\omega _0\in \mathscr P(\rr {2d})$ and invariant
BF-space $\mathscr B$.

\vrum

\item If $a\in S_{(\omega )}(\rr {2d})$ and $\op _t(a)$ is an 
isomorphism from $M^2_{(\omega _1)}(\rr d)$ to $M^2_{(\omega
_1/\omega )} (\rr d)$ for some $\omega _1\in \mathscr P(\rr
{2d})$, then $\op _t(a)$ is an isomorphism from $M_{(\omega
_2)}(\mathscr B)$ to $M_{(\omega _2/\omega )} (\mathscr B)$, for
every  $\omega _2\in \mathscr P(\rr {2d})$ and every  invariant
BF-space $\mathscr B$.
\end{enumerate}
\end{thm}

\par

We need some preparations for the proof. As a first step we prove that modulation
spaces of Hilbert type agree with certain types of Bony-Chemin spaces
(cf. Section \ref{sec1}).

\par

If $g$ is the (standard) Euclidean metric on $\rr {2d}$, then the definition of Bony-Chemin spaces (cf. Section \ref{sec1} and Appendix) specializes to the following condition. Let 
$ \psi \in C_0^\infty (\rr {2d})\setminus \{0\}$ be non-negative, even and supported in a ball of radius $1/4$,
and set $\tau _Y \psi =\psi
(\cdo -Y)$. Then $H(\omega ,g )=H(\omega )$ consists of all
$f\in \mathscr S'(\rr d)$ such that
\begin{equation}\label{Homeganorm}
\nm {f}{H(\omega )} = \Big (\int _{\rr {2d}}\omega (Y)^2\nm {\op
^w(\tau _Y \psi )f}{L^2}^2\, dY\Big )^{1/2}
\end{equation}
is finite.

\par

\begin{prop}\label{bonysobmod}
Assume that $\omega \in \mathscr P(\rr {2d})$. Then $H(\omega
)=M^2_{(\omega )}(\rr {2d})$ with equivalent norms.
\end{prop}

\par

We need two lemmas for the proof. Recall  that a linear
operator $T$ is of 
trace-class  if
$$
\sup \sum |(Tf_j,g_j)|<\infty ,
$$
where the supremum is taken over all orthonormal sequences $(f_j)$
and $(g_j)$ in $L^2(\rr d)$.

\par

The following result is an immediate consequence of the spectral
theorem for compact operators. We refer to \cite{Jan84} or Lemma 1.3 and
Proposition 1.10 in \cite{To3} for the proof of the
first part. The second part follows from the first part and
straight-forward computations.

\par

\begin{lemma}\label{traceclassel}
If $\psi \in \mathscr{S}(\rr {2d})$, then $\op ^w(\psi )$ is a
trace-class operator, and there exist two   orthonormal
sequences $(f_j)$ and $(g_j)$ in $\mathscr{S}(\rr 
d)$ and a non-negative non-increasing  sequence $(\lambda _j) \in \ell
^1$ such that 
$$
\op ^w(\psi )f = \sum _{j=1}^\infty \lambda _j (f, f_j) g_j \, ,
$$
when $f\in L^2(\rr d)$. Moreover, set $f_{j,Y}(x)=e^{i\scal x\eta}f_j(x-y)$ and $g_{j,Y}(x)=e^{i\scal x\eta}g_j(x-y)$ for
$Y=(y,\eta)$, then
\begin{equation}\label{eq:gh}
  \op ^w(\tau _Y \psi )f = \sum _{j=1}^\infty \lambda _j (f, f_{j,Y}) g_{j,Y} \, .
\end{equation}
\end{lemma}

\par

We also need the following lemma. Since it is difficult to find a
proof in the literature, we include its short proof.

\par

\begin{lemma}\label{schwartzproducts}
Assume that $f\in \mathscr S(\rr {d_1+d_2})$. Then there are $f_0\in
\mathscr S(\rr {d_1+d_2})$ and strictly positive rotation-invariant
function  $g\in \mathscr
S(\rr {d_1})$ such that
$$
f(x_1,x_2)=f_0(x_1,x_2)g(x_1).
$$
\end{lemma}

\par

\begin{proof}
We only prove the result for $d_1=d$ and $d_2=0$. The general case
is  similar  and  left to the reader. For each integer  $j\ge 1$
define the set  
$$
\Omega _j=\sets {x\in \rr d}{\sum _{|\alpha |,|\beta |\leq
2j}|x^\alpha D^{\beta }f(x)|\ge 2^{-2j}\eabs x^{-2j}}\, .
$$
Since $f\in \mathscr {S}(\rr {d_1})$, $\Omega _j$
is compact, and $\Omega _j \subseteq \Omega _{j+1}$ for all $j$. 

\par

Let $R_0=-1$ and
$
R_j=j+\sup \sets {|x|}{x\in  \Omega _j}$ for $ j\ge 1,
$ and let $(\fy _j)_{j=0}^\infty$ be a bounded set in
$C_0^{\infty}(\mathbf R)$ such that $\fy _j\ge 0$,
\begin{align*}
\supp \fy _j  &\subseteq \sets {r}{R_j-1 \le r\le R_{j+1}+1}
\intertext{and}
\sum _{j=0}^\infty \fy _j(r) &= 1 \quad \text{when}\quad
r\ge 0.
\end{align*}
Now set 
$$
g(x)=\sum _{j=0}^\infty \fy _j(|x|)2^{-j}\eabs x^{-j},\quad
\text{and}\quad f_0(x)=f(x)/g(x).
$$
Then $f_0, g\in \mathscr{S}(\rr {d})$ and $f= f_0 g$. 
\end{proof}

\par

\begin{proof}[Proof of Proposition \ref{bonysobmod}.]
Let $\psi \in C_0^\infty (\rr {2d})\setminus \{0\}$ be  as in
\eqref{Homeganorm}, and let
$$
G(x,z)=(\mathscr F_2\psi )\Big ( \frac{x+z}{2},z-x \Big ),
$$
which belongs to $\mathscr S(\rr {2d})$. By definition, $\omega$ 
is $v$-moderate for some $v\in \mathscr P(\rr {2d})$, and by
Lemma \ref{schwartzproducts} we may choose
 $G_1\in \mathscr S(\rr {2d})$ and a strictly positive 
$\fy \in \mathscr S(\rr d)$, such that
$G(x,z)=G_1(x,z)\fy (z)$. 
Using formula \eqref{eq:gr} for $f\in \mathscr S(\rr d)$ and $Y=(y,\eta
)\in \rr {2d}$, the expression for the $H(\omega )$-norm of $f$
becomes 
\begin{eqnarray*}
  \lefteqn{(2\pi )^{2d}\nm f{H(\omega )}^2} \\[1ex]
&=& \iint \Big \vert \omega (Y) \iint \psi \Big (\frac {x+z}2-y,\xi
-\eta \Big )f(z)e^{i\scal {x-z}\xi }\, dzd\xi \Big \vert ^2\, dxdY
\\[1ex]
&=& \iint \Big \vert \omega (Y) \iint \psi \Big (\frac {x+z}2,\xi
\Big )f(z+y)e^{-i\scal z\eta} e^{i\scal {x-z}\xi }\, dzd\xi \Big
\vert ^2\, dxdY
\\[1ex]
&=& \iint \Big \vert \omega (Y) \int G(x,z)f(z+y)e^{-i\scal z\eta}\,
dz\Big \vert ^2\, dxdY
\\[1ex]
&=& \iint \Big \vert \omega (Y) \int G_1(x,z)\fy
(z)f(z+y)e^{-i\scal z\eta}\, dz\Big \vert ^2\, dxdY.
\end{eqnarray*}
In the second equality we have taken $z-y$, $\xi -\eta$, $x-y$ and
$Y$ as the new variables of integration. Since the inner integral on
the right-hand side is the Fourier transform of the product
$G_1(x,z) \cdot \big ( \fy (z)f(z+y)\big )$ with respect to the
 variable $z$, we obtain
\begin{eqnarray*}
\lefteqn{(2\pi )^{3d}\nm f{H(\omega )}^2}
\\[1ex]
&=&  \iint \Big \vert \omega (Y) \Big(|\mathscr F_2(G_1)(x,\cdo
)|*|\mathscr F(\fy \, f(\cdo +y))|\Big)(\eta )\Big \vert ^2\, dxdY
\\[1ex]
&=& \iint \Big \vert \omega (Y) \big (|\mathscr F_2(G_1)(x,\cdo
)|*|V_\fy f(y,\cdo )|\big ) (\eta )\Big \vert ^2\, dxdY
\\[1ex]
&\le & C_1\iint \Big \vert \Big(|\mathscr F_2(G_1)(x,\cdo
)v(0,\cdo )|*|V_\fy f(y,\cdo )\omega (y,\cdo )|\Big) (\eta )\Big
\vert ^2\, dxdY
\\[1ex]
&\le & C_2\nm {V_\fy f\, \omega }{L^2}^2=C_2\nm f{M^2_{(\omega )}}^2,
  \end{eqnarray*}
for some constants $C_1$ and $C_2$. Here the last inequality
follows from Young's inequality. Hence $M^2_{(\omega )}(\rr
{2d})\subseteq H(\omega )$.

\par

For the reverse inclusion  we note that $(f_{j,Y})$ and $(g_{j,Y})$ in
the spectral representation \eqref{eq:gh} are orthonormal sequences for each fixed $Y\in \rr {2d}$. Hence \eqref{eq:gh} and Bessel's inequality give
$$
\nm {\op ^w(\tau _Y \psi )f}{L^2}  = \Big(\sum _{j=1}^\infty \lambda _j^2
|( f, f_{j,Y})_{L^2} |^2\Big)^{1/2} \ge \lambda _1 |(f, f_{1,Y})_{L^2}|=\lambda _1 |V_{f_{1}} f(y,\eta )|.
$$
A combination of these estimates gives
\begin{eqnarray*}
 \nm f{H(\omega )}^2 &=& \int \omega (Y)^2\nm {\op ^w(\tau_Y \psi
 )f}{L^2}^2\, dY \\
 & \ge &  \lambda _1^2   \iint \omega (y,\eta )^2|V_{f_1} f(y,\eta )|^2\, dyd\eta =\nm
 f{M^2_{(\omega )}}^2 \, ,   
\end{eqnarray*}
where the last identity follows from Proposition~\ref{p1.4}(1) on norm
equivalence, and the fact that $f_1$ belongs to $\mathscr S$. 
Consequently $H(\omega )\subseteq M^2_{(\omega )}(\rr {2d})$.

\par

Summing up, we have shown that  $H(\omega )= M^2_{(\omega )}(\rr {2d})$ with
equivalent norms, and the proof is complete.
\end{proof}

\par

\begin{proof}[Proof of Theorem \ref{identification}]
We first remark that by Proposition~\ref{S(omega)prop}(1) the
statements are independent of the pseudo-differential 
calculus used. 
Hence we may  assume that $t=1/2$ and use the Weyl calculus. 

\par

To prove (1), we use a fundamental result of Bony and
Chemin~\cite{BC}. According to \cite[Corollary~6.6]{BC}  and Proposition
\ref{p3.2}  there exist  $a\in
S_{(\omega )}(\rr {2d})$ and $b\in S_{(1/\omega )}(\rr {2d})$ such
that both $\op ^w(a)\circ \op ^w(b)$ and $\op ^w(b)\circ \op ^w(a)$ are
equal to the identity operator on $H(\omega _0)=M^2_{(\omega _0)}(\rr d)$
for every  $\omega _0\in
\mathscr P(\rr {2d})$. Since the symbol $a \wpr b$ of $\op
^w(a)\circ \op ^w(b)$ is in $S_{(1)}=S^0_{0,0}$ by
Proposition~\ref{S(omega)prop}(2), the boundedness result of
Proposition~\ref{p3.2} is applicable and shows that
both $\op ^w(a)\circ \op ^w(b)$ and $\op ^w(b)\circ
\op ^w(a)$ are equal to the  identity operator on $M_{(\omega _0)}(\mathscr B)$
for arbitrary  $\omega _0\in \mathscr P(\rr {2d})$
and invariant BF-space $\mathscr B$.

\par

Since $\op ^w(a)$ maps $M_{(\omega _0)}(\mathscr B)$ to
$M_{(\omega _0/\omega )}(\mathscr B)$ and $\op ^w(b)$ maps
$M_{(\omega _0/\omega)}(\mathscr B)$ to $M_{(\omega _0 )}(\mathscr B)$, the factorization
of the identity operator implies that these mappings are one-to-one and
onto. Hence $ \op ^w(a)$ is an isomorphism
between $M_{(\omega _0)}(\mathscr B)$ and $M_{(\omega _0/\omega )}(\mathscr B)$.

\par

(2) By (1), we may find
$$
a_1\in S_{(\omega _1)},\quad b_1\in S_{(1/\omega _1)},\quad a_2\in
S_{(\omega _1/\omega )},\quad b_2\in S_{(\omega /\omega _1)}
$$
satisfying  the following properties:
\begin{itemize}
\item $\op ^w(a_j)$ and $\op ^w(b_j)$ are inverses to each others on
$\mathscr S'(\rr d)$ for $j=1,2$;

\vrum

\item For arbitrary  $\omega _2\in \mathscr P(\rr {2d})$, each of  the mappings
\begin{equation}\label{4homeomorphisms}
\begin{aligned}
\op ^w(a_1)\, &:\, M^2_{(\omega _2)}\to M^2_{(\omega
_2/\omega _1)},
\\[1ex]
\op ^w(b_1)\, &:\, M^2_{(\omega _2)}\to M^2_{(\omega
_2\omega _1)},
\\[1ex]
\op ^w(a_2)\, &:\, M^2_{(\omega _2)}\to M^2_{(\omega
_2\omega /\omega _1)},
\\[1ex]
\op ^w(b_2)\, &:\, M^2_{(\omega _2)}\to M^2_{(\omega
_2\omega _1/\omega )}
\end{aligned}
\end{equation}
is an  isomorphism.
\end{itemize}

\par

In particular, $\op ^w(a_1)$ is an isomorphism from $M^2_{(\omega _1 )}$ to
$L^2$, and $\op ^w(b_1)$ is an isomorphism  from $L^2$ to
$M^2_{(\omega _1)}$. 

\par

Now set $c=a_2\wpr a \wpr b_1$. Then by Proposition \ref{S(omega)prop}(2), the symbol $c$
satisfies  
$$
c=a_2\wpr a \wpr b_1\in S_{(\omega _1/\omega )}\wpr S_{(\omega )}\wpr
S_{(1/\omega _1)}\subseteq S_{(1)}=S_{0,0}^0. 
$$
Furthermore,  $\op ^w(c)$ is a composition of three isomorphisms and
consequently  $\op ^w(c)$ is boundedly invertible on  $L^2$. 
\par 
By the
Wiener algebra property of $S_{0,0}^0$ with respect to the Weyl product
(cf. Proposition~\ref{specinv} or  \cite{Beal,GR08}), the inverse of $\op ^w(c)$
is equal to $\op ^w(c_1)$ for some $c_1\in S_{0,0}^0$. Hence, by (1) it follows that $\op ^w(c)$
and $\op ^w(c_1)$ are isomorphisms on $M_{(\omega _2)}(\mathscr
B)$, for each $\omega _2\in \mathscr P(\rr {2d})$. Since  $\op ^w(c)$
and  $\op ^w(c_1)$ are bounded on every $M_{(\omega)}(\mathscr{B})$,
the factorization of the identity $ \op ^w(c) \op ^w(c_1) = \operatorname{Id}$
is well-defined on every $M_{(\omega)}(\mathscr{B})$. Consequently, $
\op ^w(c)$ is an isomorphism on $M_{(\omega)}(\mathscr{B})$.

Using the inverses of $a_2$ and $b_1$, we now find that 
$$
\op ^w(a)=\op ^w(a_1)\circ \op ^w(c)\circ \op ^w(b_2) 
$$
is a composition of  isomorphisms   from the domain space  $M_{(\omega
  _2)}(\mathscr B)$ onto the target space 
$M_{(\omega _2/\omega )} (\mathscr B)$ (factoring through some
intermediate spaces)  for every  $\omega _2\in
\mathscr P(\rr {2d})$ and every translation invariant BF-space $\mathscr
B$. The proof is complete.
\end{proof}

\par

\section{Mapping properties for Toeplitz operators}\label{sec2}

\par

In this section we study the  bijection properties of  Toeplitz
operators between modulation spaces. We  first
state  results for  Toeplitz operators that  are well-defined in
the sense of \eqref{toeplitz} and Propositions \ref{Tpcont1} and
\ref{Tpcont2}. Then we state and prove more general results
for  Toeplitz operators that  are   defined  only  in the framework of
pseudo-differential calculus.

\par

We start with the following result about  Toeplitz operators with smooth symbols.

\par

\begin{thm}\label{locidentification}
Assume that $\omega ,v\in \mathscr P(\rr {2d})$, $\omega _0\in
\mathscr P_0(\rr {2d})$, and  that
 $\omega _0$ is $v$-moderate. If $\fy \in M^1_{(v)}(\rr d)$ and 
$\mathscr B$ is a translation invariant BF-space,  then $\tp _\fy
(\omega _0)$ is an isomorphism from $M_{(\omega )}(\mathscr B)$ to
$M_{(\omega /\omega _0)}(\mathscr B)$.
\end{thm}

\par

In the next result we relax our restrictions on the weights but impose
more restrictions on the modulation spaces.

\par

\begin{thm}\label{locidentification2}
Assume that $0\le t\le 1$, $p,q\in [1,\infty]$, and $\omega ,\omega
_0,v_0,v_1\in \mathscr P(\rr {2d})$ are such that  $\omega _0$ is $v_0$-moderate and  $\omega $ is
$v_1$-moderate. Set  $v=v_1^tv_0$, $\vartheta =\omega _0^{1/2}$ and let $\omega _{0,t}$ be the same as in \eqref{omega0t}.
If  $\fy \in M^{1}_{(v )}(\rr d)$ and $\omega _0\in
\mathcal M^{\infty}_{(1/\omega _{0,t})}(\rr {2d})$,  then $\tp _\fy (\omega _0)$ is an
isomorphism from $M_{(\vartheta
\omega )}^{p,q}(\rr d)$ to $M_{(\omega /\vartheta )}^{p,q}(\rr
d)$.
\end{thm}

\par

Before the proofs of Theorems \ref{locidentification} and
\ref{locidentification2} we state  the following consequence of Theorem
\ref{locidentification2} which  was the original  goal of our
investigations.

\par

\begin{cor}\label{locidentification3}
Assume that $\omega ,\omega _0,v_1,v_0\in \mathscr P(\rr {2d})$ and
that   $\omega _0$ is
$v_0$-moderate and $\omega $ is $v_1$-moderate. Set  $v=v_1v_0$
and $\vartheta =\omega _0^{1/2}$. If  $\fy \in
M^{1}_{(v)}(\rr d)$, then  $\tp
_\fy (\omega _0)$ is an isomorphism from $M_{(\vartheta \omega
)}^{p,q}(\rr d)$ to $M_{(\omega /\vartheta )}^{p,q}(\rr d)$
simultaneously for all  $p,q\in [1,\infty ]$.
\end{cor}

\par

\begin{proof}
Let $\omega _1\in \mathscr P_0(\rr {2d})$ be such that
$C^{-1}\le\omega _1/\omega _0\le C$, for some constant $C$. Hence,
$\omega /\omega _0\in L^\infty \subseteq M^{\infty}$. By Theorem 2.2
in \cite{To9}, it follows that $\omega =\omega _1 \cdot (\omega
/\omega _1)$ belongs to $M^{\infty}_{(\omega _2)}(\rr {2d})$, when
$\omega _2(x,\xi ,\eta ,y)=1/\omega _0(x,\xi )$. The result now
follows by setting $t=1$ and $q_0=1$ in Theorem
\ref{locidentification2}.
\end{proof}

\par

In the proofs of Theorems \ref{locidentification} and
\ref{locidentification2} we consider Toeplitz operators as defined by
an extension of the form~\eqref{toeplitz}.
Later on we  present extensions of
these theorems (cf. Theorems \ref{locidentification}$'$ and
\ref{locidentification2}$'$ below) for those readers who accept to use pseudo-differential
calculus to extend the definition of Toeplitz operators. For the proofs of Theorems \ref{locidentification} and
\ref{locidentification2} we therefore refer to the proofs of these extensions.

\par

We need some preparations and start with the following lemma.

\par

\begin{lemma}\label{Abijections}
Let  $\omega ,v\in \mathscr P(\rr {2d})$ be such  that $\vartheta =
\omega ^{1/2}$ is $v$-moderate.  Assume that $\fy \in M^2_{(v)}$.
Then $\tp _\fy (\omega )$  is an isomorphism from
$M^2_{(\vartheta )}(\rr d)$ onto $M^2_{(1/\vartheta )}(\rr d)$.
\end{lemma}

\par

\begin{proof}
Recall from Remark \ref{extwindows}  that for  $\fy \in M^2_{(v)}(\rr
d)\setminus \{0\}$ the expression $\nm {V_\fy f\cdot \vartheta }{L^2}$
defines an  equivalent norm
on $M^2_{(\vartheta )}$. Thus the occurring  STFTs with respect to
$\fy  $ are well defined. 

\par

Since $ \tp _\fy (\omega ) $ is bounded from  $M^2_{(\vartheta )}$ to
$M^2_{(1/\vartheta )}$ by Proposition~\ref{Tpcont2}, the estimate
\begin{equation}
  \label{eq:23}
\|\tp _\fy (\omega )  f \|_{M^2_{(1/\vartheta )}} \leq C \|f\|_{M^2_{(\vartheta )}}\,  
\end{equation}
holds for all $f\in M^2_{(\vartheta )}$. 

\par

Next we  observe that
\begin{equation}\label{Aomegaident}
(\tp _\fy (\omega ) f,g)_{L^2(\rr d)} = (\omega V_\fy f,V_\fy
g)_{L^2(\rr {2d})} = (f,g)_{M^{2,\fy }_{(\vartheta )}},
\end{equation}
for  $f,g\in M^{2}_{(\vartheta )}(\rr d)$ and $\fy \in M^2_{(v)}(\rr
d)$. The duality of modulation spaces (Proposition~\ref{p1.4}(3)) now yields the following identity:
\begin{eqnarray}
  \|f\|_{M^{2}_{(\vartheta )}} &=& \sup _{\nm g{M^{2}_{(\vartheta )}}=1}
    |(f,g)_{M^2_{(\vartheta)}}| \notag \\
&=& \sup _{\nm g{M^{2}_{(\vartheta )}}=1} |(\topo f, g )_{L^2}|
= \|\topo f\|_{M^2_{(1/\vartheta )}} \label{eq:24} \, .
\end{eqnarray}
A combination of \eqref{eq:23} and \eqref{eq:24} shows that $\nm
f{M^2_{(\vartheta )}}$ and $\| \topo 
\|_{M^2_{(1/\vartheta )}}$ are equivalent norms on $M^2_{(\vartheta )}$. 

\par

In particular, $\tp _\fy (\omega )$ is one-to-one on 
$M^2_{(\vartheta )}$ and has closed range in  $M^2_{(1/\vartheta )}$.
Since $\tp _\fy (\omega )$ is self-adjoint with respect to $L^2$, it
follows by duality that  $\tp _\fy (\omega )$ has dense
range in $M^2_{(1/\vartheta )}$. Consequently, $\topo $ is onto
$M^2_{(1/\vartheta )}$. By Banach's theorem $\topo $ is an isomorphism from
$M^2_{(\vartheta )}$ to $M^2_{(1/\vartheta )}$. 
\end{proof}

\par

We  need a further  generalization  of Proposition
\ref{Tpcont1} to more general classes of symbols and windows. Set 
\begin{equation}\label{Tomega}
\omega _1(X,Y)=\frac{v_0(2Y)^{1/2}v_1(2Y)}{\omega
_0(X+Y)^{1/2}\omega _0(X-Y)^{1/2}}.
\end{equation}

\par

\renewcommand{\rubrik}{Proposition \ref{Tpcont1}$'$}

\begin{tom}
Let $0\le t\le 1$, $p,q,q_0\in [1,\infty]$, and  $\omega ,\omega _0,v_0,v_1\in
\mathscr P(\rr {2d})$ be such that $v_0$ and $v_1$ are submultiplicative, $\omega _0$ is
$v_0$-moderate and $\omega $ is $v_1$-moderate. Set 
\begin{equation*}
r_0=2q_0/(2q_0-1),\quad v=v_1^tv_0 \quad \text{and}\quad \vartheta = \omega _0^{1/2} \,,  
\end{equation*}
and  let  $\omega _{0,t}$ and $\omega _1$ be as in \eqref{omega0t} and
\eqref{Tomega}. Then the following are true:

\par

\begin{enumerate}
\item The definition of $(a,\fy )\mapsto \tp _\fy (a)$ from $\mathscr
S(\rr {2d})\times \mathscr S(\rr d)$ to $\mathcal L(\mathscr S(\rr
d),\mathscr S'(\rr d))$ extends uniquely to a continuous map from
$\mathcal M^{\infty ,q_0}_{(1/\omega _{0,t})}(\rr {2d})\times
M^{r_0}_{(v)}(\rr d)$ to $\mathcal L(\mathscr S(\rr d),\mathscr
S'(\rr d))$.

\vrum

\item If $\fy \in M^{r_0}_{(v)}(\rr d)$ and $a\in \mathcal
M^{\infty ,q_0}_{(1/\omega _{0,t})}(\rr {2d})$, then $\tp _\fy (a)
=\op ^w(a_0)$ for some $a_0\in \mathcal M^{\infty ,1}_{(\omega _1)}(\rr
{2d})$, and $\tp _\fy (a)$ extends uniquely to a continuous map
from $M_{(\vartheta \omega )}^{p,q}(\rr d)$ to $M_{(\omega
/\vartheta )}^{p,q}(\rr d)$.
\end{enumerate}
\end{tom}

\par

For the proof we need the following special case of \cite[Proposition 2.1]{To8}.

\par

\begin{lemma}\label{Prop2.1inTo8}
Assume that $q_0,r_0\in [1,\infty ]$ satisfy $r_0=2q_0/(2q_0-1)$. Also
assume that $v \in \mathscr P(\rr {2d})$ is submultiplicative, and
that $\kappa ,\kappa _0\in \mathscr P(\rr {2d}\oplus \rr {2d})$
satisfy 
\begin{equation}\label{kappacond1}
\kappa _0(X_1+X_2,Y)\le C\kappa (X_1,Y)\, v(Y+X_2)v(Y-X_2)\quad X_1,X_2,Y\in \rr {2d},
\end{equation}
for some constant $C>0$. Then the  map $(a,\fy )\mapsto \operatorname 
{Tp}_{\fy} (a)$ from $\mathscr S(\rr {2d})\times \mathscr S(\rr {d})$
to $\mathcal L(\mathscr S(\rr {d}),\mathscr S'(\rr {d}))$ extends
uniquely to a continuous mapping from $\mathcal M^{\infty
  ,q_0}_{(\omega )}(\rr {2d})\times M^{r_0}_{(v)}$  to $\mathcal L(\mathscr S(\rr
{d}),\mathscr S'(\rr {d}))$. Furthermore, if $\fy \in
M^{r_0}_{(v)}(\rr d)$ and $a\in \mathcal M^{\infty ,q_0}_{(\kappa 
)}(\rr {2d})$, then $\operatorname {Tp}_{\fy} (a) = \op ^w(b)$ for some $b\in \mathcal M^{\infty ,1}_{(\kappa _0)}$.
\end{lemma}

\par

\begin{proof}[Proof of Proposition \ref{Tpcont1}{\,}${}^\prime$]
We show that the conditions on the involved parameters and weight
functions satisfy  the conditions of  Lemma \ref{Prop2.1inTo8}.

\par

First we observe that
$$
v_j(2Y)\le Cv_j(Y+X_2)v_j(Y-X_2),\quad j=0,1
$$
for some constant $C$ which is independent of $X_2,Y\in \rr {2d}$,
because $v_0$ and $v_1$ are submultiplicative. Refering back to \eqref{Tomega} this gives
\begin{multline*}
\omega _1(X_1+X_2,Y) = \frac {v_0(2Y)^{1/2}v_1(2Y)}{\omega
_0(X_1+X_2+Y)^{1/2}\omega _0(X_1+X_2-Y)^{1/2}}
\\[1ex]
\le C_1 \frac
{v_0(2Y)^{1/2}v_1(2Y)v_0(X_2+Y)^{1/2}v_0(X_2-Y)^{1/2}}{\omega
_0(X_1)}
\\[1ex]
= C_1 v_1(2Y)^{1-t}\frac
{v_0(2Y)^{1/2}v_1(2Y)^t v_0(X_2+Y)^{1/2}v_0(X_2-Y)^{1/2}}{\omega
_0(X_1)}
\\[1ex]
\le C_2v_1(2Y)^{1-t} \frac {v_1(X_2+Y)^t v_1(X_2-Y)^t v_0(X_2+Y)v_0(X_2-Y)}{\omega
_0(X_1)}.
\end{multline*}
Hence
\begin{equation}\label{omegacond1}
\omega _1(X_1+X_2,Y) \le C\frac
{v_1(2Y)^{1-t}v(X_2+Y)v(X_2-Y)}{\omega _0(X_1)}.
\end{equation}
By letting $\kappa _0 = \omega _1$ and $\kappa =1/\omega _{0,t}$, it
follows that \eqref{omegacond1} agrees with \eqref{kappacond1}. The
result now follows from Lemma \ref{Prop2.1inTo8}. 
\end{proof}

\par

In the remaining part  of the paper we interpret  $\tp
_\fy (a)$  as the extension  of a  Toeplitz
operator provided by Proposition \ref{Tpcont1}$'$. (See also Remark
\ref{extensionremark} below for more comments.)

\par

Proposition \ref{Tpcont1}$'$ can also be applied on Toeplitz operators with smooth weight.

\par

\begin{prop}\label{Aomegaproperties}
Assume that $\omega _0\in \mathscr P_0(\rr {2d})$, that $v\in \mathscr
P(\rr {2d})$ is submultiplicative,  and that   $\omega _0
^{1/2}$ is $v$-moderate. If $\fy \in M^2_{(v)}$, then   $\tp _\fy
(\omega _0) =\op ^w(b)$ for some $b\in S_{(\omega _0
)}(\rr {2d})$.
\end{prop}

\par

\begin{proof}
By Proposition~\ref{S(omega)prop}  we have $\omega _0 \in \mathcal
M^{\infty ,1}_{(1/\omega _{0,N})}(\rr {2d})$ for every  $N\ge 0$, where
$\omega _{0,N}(X,Y)=\omega _0(X)\eabs Y^{-N}$. Furthermore, 
$$
\omega _1(X,Y)=\frac{\eabs Y^{N}v(2Y)^{1/2}}{\omega
_0(X+Y)^{1/2}\omega _0(X-Y)^{1/2}} \ge C\eabs Y^{N-N_0}/\omega _0(X),
$$
for some constants $C$ and $N_0$ which are independent of $N$.
Proposition \ref{Tpcont1}$'$ implies that existence of  some $b \in
\mathcal M^{\infty ,1}_{(1/\omega _{0,N})}$, such
that   $\tp _\fy (\omega _0) =\op ^w(a)$. Applying  Proposition
\ref{S(omega)prop} (3) once again, we find that $b\in \bigcap _{N\geq
  0} \mathcal M^{\infty ,1}_{(1/\omega _{0,N})} = S_{(\omega
  )}(\rr{2d})$. 
\end{proof}

\par

\begin{rem}\label{extensionremark}
As remarked and stated before, there are different ways to extend the
definition of  a Toeplitz operator $\tp _\fy (a)$  (from  $\fy \in \mathscr
S(\rr d)$ and $a\in \mathscr S(\rr {2d})$) to more general classes of
symbols and windows. For example, Propositions
\ref{Tpcont1} and \ref{Tpcont2} are  based on the ``classical''
definition \eqref{toeplitz} of such operators and  a straight-forward
extension  of the $L^2$-form on $\mathscr
S$. Proposition~\ref{Aomegaproperties} interprets  $\topo $ as a
\psdo .  Let us  emphasize
that in this  context the bilinear form \eqref{toeplitz} may not be well
defined,  even  when $\fy \in M^2_{(v)}(\rr d)$ and  $\omega \in \mathscr P_0(\rr {2d})$.

\par

To shed some light on this subtlety, consider a window $\fy \in L^2
\setminus M^1$ with normalization $\|\fy \|_{L^2} =1 $ and the symbol
$\omega  \equiv 1$. Then the corresponding Toeplitz operator   $\tp
_\fy (\omega )$ is the identity operator. This is nothing but the
inversion formula for the short-time Fourier transform, e.g.,~\cite{Gc2}. Clearly the
identity operator is an isomorphism on every space. However, the
Toeplitz operator in \eqref{toeplitz}, $\tp _\fy (\omega )$ is not
defined on $M^\infty $ because it is not clear what $( 1\cdot V_\fy f
, V_\fy g)$ from \eqref{toeplitz} means for $\fy \in L^2$,
$f\in M^\infty$ and $g \in M^1$.

\par

In Theorems \ref{locidentification}$'$ and \ref{locidentification2}$'$
below, we will  extend the definition of Toeplitz operators within the
framework of pseudo-differential calculus and  we 
interpret  Toeplitz operators as pseudo-differential operators. With
this understanding, the stated mapping properties are well-defined.

\par

The reader who is not interested in full generality  or does not accept Toeplitz
operators that  are not defined directly by an extension of
 \eqref{toeplitz} may only consider the case when the
windows belong to $M^1_{(v)}$. 
For the more general window classes in Theorems
\ref{locidentification}$'$ and \ref{locidentification2}$'$ below,
however,  one
should then interpret the involved operators as ``pseudo-differential
operators that extend Toeplitz operators''.
\end{rem}

\par

The following generalization of Theorem \ref{locidentification} is
an immediate consequence of Theorem \ref{identification} and Proposition
\ref{Aomegaproperties}.

\par

\renewcommand{\rubrik}{Theorem \ref{locidentification}$'$}

\begin{tom}
Assume that $\omega ,v\in \mathscr P(\rr {2d})$, $\omega _0\in \mathscr
P_0(\rr {2d})$ and  that $\omega _0$ is $v$-moderate.
If $\fy \in M^2_{(v)}(\rr d)$ and $\mathscr B$ is
a translation invariant BF-space, then $\tp _\fy (\omega _0)$ is an
isomorphism from $M_{(\omega )}(\mathscr B)$ to $M_{(\omega /\omega
_0)}(\mathscr B)$.
\end{tom}

\par

Theorem \ref{locidentification}$'$ holds only  for smooth weight functions. In order to   relax
the conditions on the weight function
$\omega _0$,  we use the  Wiener algebra property of  $\mathcal M^{\infty ,1}_{(v)}$ instead of
$S^0_{0,0}$. On the other hand, we have  to restrict our results
to modulation spaces of the form $M^{p,q}_{(\omega )}$ instead of
$M_{(\omega )}(\mathscr{B})$.

\par

\renewcommand{\rubrik}{Theorem \ref{locidentification2}$'$}

\begin{tom}
Assume that $0\le t\le 1$, $p,q,q_0\in [1,\infty]$ and $\omega ,\omega
_0,v_0,v_1\in \mathscr P(\rr {2d})$ are such that $\omega _0$ is
$v_0$-moderate and  $\omega $ is
$v_1$-moderate. Set  $r_0=2q_0/(2q_0-1)$,  $v=v_1^tv_0$, $\vartheta =\omega _0^{1/2}$ and
let $\omega _{0,t}$ be the same as in \eqref{omega0t}.
If  $\fy \in M^{r_0}_{(v )}(\rr d)$ and $\omega _0\in \mathcal
M^{\infty ,q_0}_{(1/\omega _{0,t})}$,  then $\tp _\fy (\omega _0)$ is
an isomorphism from $M_{(\vartheta
\omega )}^{p,q}(\rr d)$ to $M_{(\omega /\vartheta )}^{p,q}(\rr d)$.
\end{tom}

\par

\begin{proof}
First we note that the Toeplitz
operator $\tp _\fy (\omega _0)$ is an isomorphism from
$M^2_{(\vartheta )}$ to $M^2_{(1/\vartheta )}$ in view of Lemma
\ref{Abijections}. With $\omega
_1$ defined in~\eqref{Tomega}, Proposition \ref{Tpcont1}$'$ implies
that there  exist  $b\in
\mathcal M^{\infty ,1}_{(\omega 
_1)}$ and $d\in \mathscr S'(\rr {2d})$ such that
$$
\tp _\fy (\omega _0) = \op ^w(b)\quad \text{and}\quad \tp _\fy (\omega
_0)^{-1}=\op ^w(d) \, .
$$
Let $\omega _2$ be the ``dual'' weight defined as 
\begin{equation}
  \label{eq:dmod1'}
  \omega _2(X,Y) =\vartheta (X-Y)\vartheta (X+Y) v_1(2Y).
\end{equation}
We will prove that $d\in \mathcal M^{\infty ,1}_{(\omega _2)}(\rr
{2d})$. 
Let us assume for now that we have already proved the existence of
such a symbol $d$. Then  we may proceed as follows. 

After checking~\eqref{e5.9}, we can apply  Proposition \ref{pseudomod} and
find that each of  the mappings
\begin{equation}\label{op(b)(d)}
\op ^w(b)\, :\, M^{p,q}_{(\omega \vartheta )} \to M^{p,q}_{(\omega /\vartheta )}\quad \text{and}\quad
\op ^w(d)\, :\, M^{p,q}_{(\omega /\vartheta )} \to M^{p,q}_{(\omega \vartheta )}
\end{equation}
is  well-defined and continuous.

\par

In order to apply Proposition~\ref{Weylprodmod}, we next check
condition~\eqref{weightprodmod}  for the weights $\omega _1$,  $\omega _2$, and
$$
\omega _3(X,Y) = \frac {\vartheta (X+Y)}{\vartheta (X-Y)}.
$$
In fact, for some constant $C_1>0$  we have
\begin{eqnarray*}
\lefteqn{\omega _1(X-Y+Z,Z)\omega _2(X+Z,Y-Z)}\\[1ex]
&=& \Big (\frac {v_0(2Z)^{1/2}v_1(2Z)}{\vartheta (X-Y+2Z)\vartheta
  (X-Y)}\Big )\cdot \big (\vartheta (X-Y+2Z)\vartheta (X+Y)
v_1(2(Y-Z)) \big ) 
\\[1ex]
&= &  \frac {v_0(2Z)^{1/2}v_1(2Z)v_1(2(Y-Z))\,  \vartheta (X+Y)}{\vartheta (X-Y)}
\\[1ex]
&\ge & C_2 \frac {\vartheta (X+Y)}{\vartheta (X-Y)} =C_2\omega _3(X,Y)\, .
\end{eqnarray*}
Therefore  Proposition \ref{Weylprodmod}  shows that the Weyl symbol
of $\op ^w(b)\circ \op ^w(d)$ belongs to $\mathcal M^{\infty
  ,1}_{(\omega _3)}(\rr {2d})$, or equivalently, $b\wpr d\in \mathcal
M^{\infty ,1}_{(\omega _3)}$. Since $\op ^w(b)$ is an isomorphism from
$M^2_{(\vartheta )}$ to $M^2_{(1/\vartheta )}$ with inverse $\op
^w(d)$, it follows that $b\wpr d =1$ and that the constant symbol $1$  belongs to $\mathcal
M^{\infty ,1}_{(\omega _3)}$. By similar arguments it follows that
$d\wpr b=1$. Therefore the identity operator $\mathrm{Id}= \op ^w (b)
\op ^w(d) $ on $M^{p,q}_{(\omega \vartheta )}$ factors through
$M^{p,q}_{(\omega /\vartheta )}$,  and thus $\op ^w(b)= \tp _\fy (\omega
_0)$ is an isomorphism
from $M^{p,q}_{(\omega \vartheta )}$ onto $M^{p,q}_{(\omega /\vartheta )}$
with inverse $\op ^w(d)$. This proves  the assertion. 

\par

It remains to prove that $d\in \mathcal M^{\infty ,1}_{(\omega
 _2)}(\rr {2d})$. Using once again the basic result of Bony and
Chemin~\cite{BC}, we  choose $a\in S_{(1/\vartheta )}(\rr {2d})$ and
$c\in 
S_{(\vartheta )}(\rr {2d})$ such that the map
$$
\op ^w(a)\, :\, L^2(\rr d)\to M^2_{(\vartheta )}(\rr d)
$$
is an isomorphism with inverse $\op ^w(c)$. By Proposition
\ref{S(omega)prop}  $\op ^w(a)$ is
also bijective from $M^2_{(1/\vartheta )}(\rr d)$ to $L^2(\rr
d)$. Furthermore, by Theorem~\ref{identification} it follows that $a\in \mathcal
M^{\infty ,1}_{(\vartheta _N)}$ for each $N\ge 0$,  
where
$$
\vartheta _N(X,Y)=\vartheta (X)\eabs Y^N.
$$
Let $b_0=a\wpr b \wpr a$. From  Proposition~\ref{corweyl} we know that 
\begin{equation}\label{b0mod}
b_0\in \mathcal M^{\infty ,1}_{(v_2)}(\rr {2d}),\quad \text{where}\quad v_2(X,Y)=v_1(2Y)
\end{equation}
is submultiplicative and depends on $Y$ only. 
Since $\op ^w(b)$ is
bijective from $M^2_{(\vartheta )}$ to $M^2_{(1/\vartheta )}$ by
Lemma \ref{Abijections} (2),  $\op ^w(b_0)$ is  bijective and
continuous  on  $L^2$. 

\medspace

Since $v_2$ is submultiplicative and in $\mathscr{P}(\rr{2d})$,  $\mathcal M^{\infty
  ,1}_{(v_2)}$ is a 
Wiener algebra by  Theorem \ref{specinv}.  Therefore, the bijective
operator $\op ^w(b_0)$   on $L^2$ possesses an inverse   $\op ^w(d_0)$
for some $d_0\in \mathcal M^{\infty 
,1}_{(v_2)}(\rr {2d})$. 

\par

Since
$$
\op ^w (d_0) = \op ^w (b_0)\inv = \op ^w(a)\inv \op ^2(b)\inv
\op ^w(a)\inv ,
$$
we find that 
$$
\op ^w (d) = \op ^w (b)\inv = \op ^w(a) \op ^w(d_0)\op ^w(a),
$$
or equivalently, 
\begin{equation}\label{d0dcond}
d = a\wpr d_0\wpr a ,\quad \text{where} \quad a\in S_{(1/\vartheta )}
\,\, \text{and}\ d_0\in \mathcal M^{\infty ,1}_{(v_2)}\, .
\end{equation}
The definitions of the weights are chosen such that  Proposition 
\ref{corweyl} implies that $d \in  \mathcal M^{\infty ,1}_{(\omega _4)}$. With
this fact, the proof is now complete.
\end{proof}

\par

\section{Examples on bijective pseudo-differential
operators on modulation spaces}\label{sec3}

\par

In this section we construct explicit isomorphisms between modulation
spaces with different weights. Applying the  results of the previous
sections, these may be in the form of  pseudo-differential operators
or of Toeplitz operators.

\par

In fact, the following two propositions are
immediate consequences of  \eqref{toeplweyl}, and Theorems
\ref{locidentification}$'$ and \ref{locidentification2}$'$. 

\par

\begin{prop}\label{isompseudos}
Assume that $\omega ,v\in \mathscr P(\rr {2d})$, $\omega _0\in \mathscr
P_0(\rr {2d})$ and  that $\omega _0$ is $v$-moderate.
If $\fy \in M^2_{(v)}(\rr d)$ and $\mathscr B$ is
a translation invariant BF-space, then $\tp _\fy (\omega _0)$ is an
isomorphism from $M_{(\omega )}(\mathscr B)$ to $M_{(\omega /\omega
_0)}(\mathscr B)$.
\end{prop}

\par

\begin{prop}\label{isompseudos2}
Assume that $0\le t\le 1$, $p,q,q_0\in [1,\infty]$ and $\omega ,\omega
_0,v_0,v_1\in \mathscr P(\rr {2d})$ are such that $\omega _0$ is
$v_0$-moderate and  $\omega $ is
$v_1$-moderate. Set  $r_0=2q_0/(2q_0-1)$,  $v=v_1^tv_0$ and
let $\omega _{0,t}$ be the same as in \eqref{omega0t}.
If  $\fy \in M^{r_0}_{(v )}(\rr d)$ and $\omega _0\in \mathcal
M^{\infty ,q_0}_{(1/\omega _{0,t})}$,  then $\op ^w(W_{\fy ,\fy}  *\omega _0)$ is
an isomorphism from $M_{(\omega _0^{1/2}
\omega )}^{p,q}(\rr d)$ to $M_{(\omega /\omega _0^{1/2} )}^{p,q}(\rr d)$.
\end{prop}

\par

\begin{cor}\label{isompseudos3}
Assume that $p,q\in [1,\infty ]$, $\omega _0,\omega \in \mathscr P(\rr
{2d})$, and let $\mathscr B$ be a translation invariant BF-space on
$\rr {2d}$. For $ \lambda = (\lambda _1$, $\lambda _2) \in \rr{2}_+ $  let $\Phi _{\lambda}  $ be the Gaussian 
$$
\Phi _\lambda  (x,\xi ) = Ce^{-(\lambda _1|x|^2+\lambda _2|\xi |^2)}\, .
$$

\par

\begin{enumerate}
\item The weight  $\omega _0 \ast \Phi _\lambda $ is  in $\mathscr
  P_0(\rr {2d})$ for all $\lambda \in \rr{2}_+ $ and
$$
C^{-1}\omega _{\Phi} \le \omega _0  \ast \Phi _\lambda \le C\omega _\Phi ,
$$
for some constant $C>0$.

\vrum

\item If $\lambda _1\cdot \lambda _2 <1$, then there exists a $\nu \in
  \rr{2}_+ $ and $\fy \in \mathscr{S}(\rr{d})$ such that $\op ^w(\omega _0 \ast \Phi _\lambda ) =
  \mathrm{Tp}_\fy (\omega _0 \ast \Phi _\nu)$ is bijective from
  $M_{(\omega )}(\mathscr B)$ to $M_{(\omega /\omega _0)}(\mathscr
  B)=M_{(\omega /\omega _\Phi)}(\mathscr B)$ for all $\omega \in \mathscr P(\rr
{2d})$.

\vrum

\item If $\lambda _1\cdot \lambda _2 \le 1$ and in addition $\omega _0
  \in \mathscr P_0(\rr {2d})$, then $\op ^w(\omega _0 \ast \Phi
  _\lambda ) = \mathrm{Tp}_\fy (\omega _0)$ is bijective from
  $M_{(\omega )}(\mathscr B)$ to $M_{(\omega /\omega _0)}(\mathscr
  B)=M_{(\omega /\omega _\Phi)}(\mathscr B)$ for all $\omega \in \mathscr P(\rr
{2d})$. 
\end{enumerate}
\end{cor}

\par

\begin{proof}
The assertion (1) follows easily from the definitions.

\par

(2) Choose $\mu _j>\lambda _j$ such that $\mu _1\cdot \mu _2=1$. Then
the Gaussian $\Phi _\mu$ is a multiple of a Wigner
  distribution, precisely $\Phi _\mu = c W(\fy , \fy )$
  with $\fy (x)=e^{-\mu _1|x|^2/2}$.   By
the semigroup property of Gaussian functions (cf. e.g., \cite{Fo,Gc2}) there
exists another Gaussian, namely $\Phi _\nu$,  such that $\Phi _\lambda
= \Phi _\mu \ast \Phi _\nu $. Using \eqref{toeplweyl}, this factorization implies that the
 Weyl operator with symbol $\omega _0 \ast \Phi _\lambda$  is in fact a Toeplitz operator, namely
 \begin{eqnarray*}
\mathrm{Op}^w (\omega _0 \ast \Phi _\lambda) &=& \mathrm{Op}^w (\omega
_0 \ast \Phi _\nu \ast \Phi _\mu) \\
& = &  \mathrm{Op}^w (\omega _0 \ast \Phi _\nu \ast cW(\fy , \fy )) \\
& = & c(2\pi )^{d/2}\mathrm{Tp}_\fy (\omega _0 \ast \Phi _\nu )   
\end{eqnarray*}

\par

By (1)  $\omega _0 \ast \Phi _\nu \in \mathscr P_0(\rr {2d})$ is
equivalent to $\omega _0$. Hence   Proposition
\ref{isompseudos} shows that 
$\mathrm{Op}^w (\omega _0 \ast \Phi _\lambda)$
is bijective from $M_{(\omega )}(\mathscr B)$ to $M_{(\omega /\omega _0)}(\mathscr B)$. This proves (2).

\par

(3) follows from (2) in the case $\lambda _1\cdot \lambda
_2<1$. If $\lambda _1\cdot \lambda _2 =1$, then as above $\Phi
_\lambda = c W(\fy , \fy )$ for $\fy(x) = e^{-\lambda_1 |x|^2/2}$ and
thus
$$
\op ^w(\omega _{0} \ast \Phi _\lambda) = \mathrm{Tp}_\fy  ^w(\omega _{0})
$$
is bijective from $M_{(\omega )}(\mathscr B)$ to $M_{(\omega /\omega
  _0)}(\mathscr B)$, since $\omega _0 \in \mathscr{P}_0(\rr{2d})$.  The proof is complete.
\end{proof}

\par

\section*{Appendix}

\par

Bony and Chemin \cite[Section 5]{BC} give a definition of the
Sobolev-type spaces  $H(\omega ,g)$ 
for a  general class of  metrics $g$ and weight functions. This  norm
is rather  complicated (cf. \cite[Section 5]{BC} for
strict definition). For example, the definition of the $H(\omega
,g)$-norm in   formula (5.1) of \cite{BC}  involves a sum of
expressions that are similar to the
right-hand side of \eqref{Homeganorm}. However, when $g$ is the
usual Euclidean metric on $\rr {2d}$, then the functions $\fy _Y$,
$\psi _{Y,\nu}$ and $\theta _{Y,\nu}$ in \cite[Definition 5.1]{BC}
can be chosen in the following way.

\par

Let $0\le \theta \in C_0^\infty (\rr {2d})\setminus 0$ be even and
supported in the ball with center at origin and radius $1/4$. Then it
follows that $\widetilde \fy =\theta *_\sigma \theta *_\sigma \theta
\in C_0^\infty (\rr {2d})\setminus 0$ is even and non-negative. Here
$*_\sigma$ is the twisted convolution, defined by the formula
$$
(a*_\sigma b)(x,\xi ) = (2/\pi )^{d/2}\iint _{\rr {2d}} a(x-y,\xi
-\eta )b(y,\eta )e^{2i (\scal y\xi -\scal x\eta )}\, dyd\eta .
$$

\par

Now let $\fy =c\widetilde \fy$, where $c>0$ is chosen such that $\nm
\fy{L^1}=1$. From Lemma 1.5 and Proposition 1.6 in \cite{To3} we have
$$
\widetilde \fy = \theta *_\sigma  \theta *_\sigma  \theta = (2\pi
)^{-d} \theta \wpr  \check \theta \wpr  \theta = (2\pi )^{-d} \theta
\wpr  \theta \wpr  \theta .
$$
By letting
\begin{alignat*}{2}
\fy _Y &= \fy (\cdo -Y),&\quad \psi _{Y,0} &= \theta _Y =\theta
(\cdo -Y),
\\[1ex]
\theta _{Y,\nu} &= \psi _{Y,\nu} =0,&\quad \nu &\ge 1,
\end{alignat*}
it follows that all the required properties in \cite[Definition
5.1]{BC} are fulfilled. Consequently, \eqref{Homeganorm} defines a
norm for $H(\omega )$.

\end{document}